\documentclass[hidelinks, 11pt]{article}
\usepackage[a4paper,bindingoffset=0.in,left=2.5cm,right=2.5cm,top=2cm,bottom=3cm]{geometry}
\usepackage{amsmath, amssymb,amsthm}
\usepackage{graphicx,subfigure}
\usepackage{color}
\usepackage{url}
\usepackage[a]{esvect}
\usepackage{enumitem}
\usepackage[utf8]{inputenc}
\usepackage{lmodern}
\usepackage{setspace}
\usepackage{wrapfig}
\usepackage{mdframed}
\usepackage[normalem]{ulem}

\newtheorem{property}{Property}
\newtheorem{proposition}{Proposition}
\newtheorem{theorem}{Theorem}

\newtheorem{definition}{Definition}

\newtheorem{corollary}{Corollary}

\usepackage{mathtools, amsmath, amssymb, graphicx, amsfonts}
\usepackage[colorlinks=true,allcolors=blue!75!black]{hyperref} %to color links
\usepackage{xcolor} %to color comments
\usepackage{bbm} %to use indicator function
\usepackage{float} %to use [h!] with figures
\usepackage{amsopn} %standard?
\usepackage{mathrsfs} %for \mathscr
\usepackage{comment} %to comment out sections
\usepackage{hyperref} %for hyperlink
\newcommand*{\colorboxed}{}
\def\colorboxed#1#{%
  \colorboxedAux{#1}%
}
\newcommand*{\colorboxedAux}[3]{%
  \begingroup
    \colorlet{cb@saved}{.}%
    \color#1{#2}%
    \boxed{%
      \color{cb@saved}%
      #3%
    }%
  \endgroup
}

\newcommand{\tr}{\operatorname{tr}}

\newcommand{\diag}{\operatorname{diag}}

\newcommand{\argmax}{\operatorname{argmax}}
\newcommand{\supp}{\operatorname{supp}}
\newcommand{\spn}{\operatorname{span}}
\newcommand{\bmu}{\boldsymbol{\mu}} %for balls 
\newcommand{\bomega}{\boldsymbol{\omega}} %for balls 
\newcommand{\bzeta}{\boldsymbol{\zeta}} %for Omega
\newcommand{\var}{\operatorname{var}}
\newcommand{\brho}{\boldsymbol{\rho}} %for balls 

\begin{document}

%\preprint{APS/123-QED}

\title{Discrete curvature on graphs from the effective resistance}

\author{
  Karel Devriendt\footnote{Contact: \href{mailto:karel.devriendt@mis.mpg.de}{karel.devriendt@mis.mpg.de} -- Funding: KD was supported by The Alan Turing Institute under the EPSRC grant EP/N510129/1.}\\
  \small{University of Oxford, UK}\\
  \small{Alan Turing Institute, UK}
  \and
  Renaud Lambiotte\footnote{Funding: RL acknowledges support from the EPSRC Grants No. EP/V013068/1 and No. EP/V03474X/1.}\\
  \small{University of Oxford, UK}\\
  \small{Alan Turing Institute, UK}
}

\date{May 2022}

%%%%%%%%%%
% ABSTRACT
%%%%%%%%%%
\maketitle
\begin{abstract}
This article introduces a new approach to discrete curvature based on the concept of effective resistances. We propose a curvature on the nodes and links of a graph and present the evidence for their interpretation as a curvature. Notably, we find a relation to a number of well-established discrete curvatures (Ollivier, Forman, combinatorial curvature) and show evidence for convergence to continuous curvature in the case of Euclidean random graphs. Being both efficient to approximate and highly amenable to theoretical analysis, these resistance curvatures have the potential to shed new light on the theory of discrete curvature and its many applications in mathematics, network science, data science and physics.
\end{abstract}

%\setcounter{tocdepth}{1}
%\tableofcontents

%%%%%%%%%%
%  BODY
%%%%%%%%%%
\section{Introduction}\label{S1: introduction}
The idea of curvature has a long and rich history in geometry and provides the mathematical language for one of the most important physical theories, Einstein's theory of general relativity, which describes the interplay between mass and energy and the (Ricci) curvature of spacetime. As illustrated in Figure \ref{fig: continuous curvature & constant curvature}, curvature is a geometric property of smooth spaces such as lines and surfaces, which aims at characterising how a space differs from flat Euclidean space and can be understood intuitively via the dispersion of its geodesics, with geodesics remaining parallel  in a flat space with zero curvature, converging in a spherical space with positive curvature or diverging in a hyperbolic space with negative curvature.
Recently, the question to formulate an analogous theory of \emph{discrete curvature} for non-smooth objects and structures such as graphs or meshes has gained a lot of attention. Some important reasons for this growing interest can be found in the ubiquity of discrete data following the advent of digital technologies as noted in \cite[Preface]{Najman_2017_modern}, in the ongoing efforts to develop a theory of quantum gravity, where questions on discrete curvature naturally arise in some attempts at combining the continuous theory of gravity and the discrete theory of quantum mechanics\footnote{There are many different approaches to quantum gravity with discreteness arising in a variety of ways; in particular, discrete curvature does not always play as important of a role as it does in \cite{Trugenberger_2017_combinatorial} or \cite{Gorard_2020_wolfram}.} \cite{Trugenberger_2017_combinatorial, Gorard_2020_wolfram}, and finally in the successful development of a number of broadly applicable notions of discrete curvature. Broadly speaking, the vision is that discrete curvature can play an important role in representing and understanding `emergent' geometric properties in discrete structures -- which at the local level appear purely combinatorial but at a larger scale feature much richer properties. Within network science for instance, important applications include community detection \cite{ni_community_2019, sia_ollivier-ricci_2019, gosztolai_unfolding_2021}, methods for graph comparison \cite{Ni_2018_alignment} and alternative ways to capture connectivity in interconnected systems \cite{weber_curvature-based_2019, murgas_quantifying_2021, jonckheere_curvature_2019, ni_ricci_2015}. In other applications, discrete curvature appears when working with finite representations of continuous or smooth structures, as in computer graphics and image processing \cite{jin_discrete_2008, Najman_2017_modern}. More recently, notions of discrete curvature have also appeared in the design and analysis of algorithms \cite{Sigbeku_2021_curved_MCMC, Topping_2021_bottlenecks_curvature,Iyer_2013_curvature_submodular, DiGiovanni_2022_curvature-aware_embedding}.

In this article, we propose a new approach to discrete curvature on graphs which differs fundamentally from existing approaches and is based on effective resistances, a concept intrinsically linked to electrical circuits and other linear systems \cite{baez_compositional_2018}, closely connected to the Laplacian matrix and with many graph-theoretic properties. The introduced \emph{resistance curvatures}, one for nodes (vertices) and one for links (edges), are very simple to state and efficient to calculate approximately, but at the same time come with a wealth of theoretical properties which make them a promising contribution to the theory of discrete curvature and of network geometry \cite{boguna_network_2021} more generally. To evidence our measure as a valid notion of curvature, we show its relation to existing notions of discrete curvature\footnote{Ollivier-Ricci and Forman-Ricci curvature are calculated with respect to additional data on the network: a metric and Markov chain for the former, and node and link weights for the latter. The relation with our proposed notion follows from a particular and somewhat atypical choice of this data, as further discussed in Section \ref{SS4.2: relation to other discrete curvatures}.} (Ollivier-Ricci, Forman-Ricci and combinatorial curvature) and we show evidence that it converges to its continuous curvature counterpart in the case of Euclidean random graphs. As an example application, we show how this curvature has a naturally associated discrete Ricci flow which is formally related to existing models of social balance in social dynamics.
\\~\\
The rest of the article is structured as follows: Section \ref{S2: preliminaries} introduces the concept of curvature and the necessary definition on graphs and the effective resistance. In Section \ref{S3: introducing resistance curvature}, the new resistance curvature is introduced with some basic properties and a number of examples. Section \ref{S4: curvature properties} discusses the different arguments that support our introduced definitions as discrete curvatures and Section \ref{S5: Discrete Ricci Flow} treats the associated discrete Ricci flow. In Section \ref{S6: conclusion}, finally, we conclude the article with a brief summary and discussion. Proofs and further technical results are presented in the Appendix.
\section{Preliminaries}\label{S2: preliminaries}
\subsection{Continuous and discrete curvature}
\begin{figure}[h!]
    \centering
    \includegraphics[width=\textwidth]{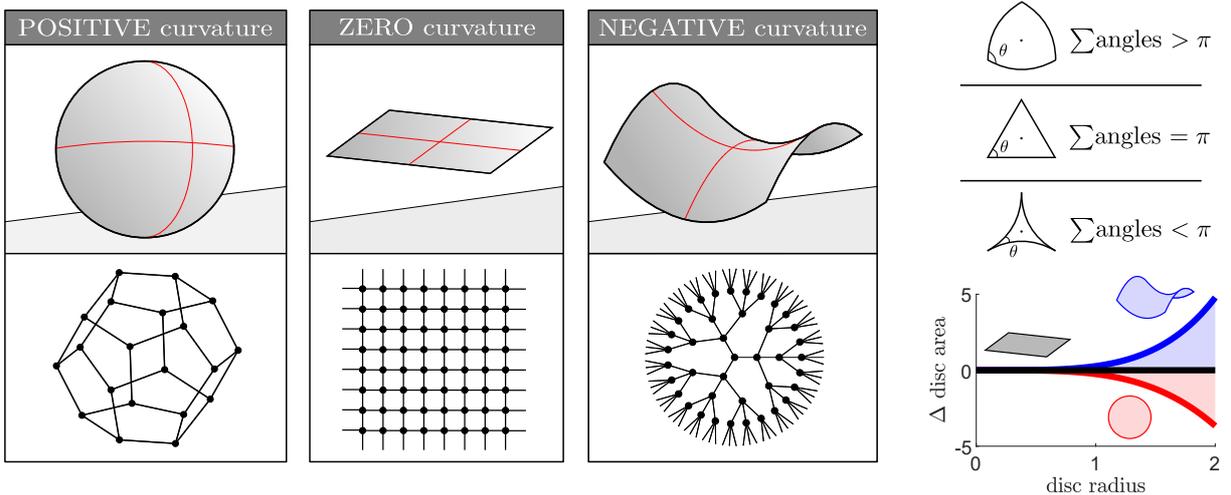}
    \caption{\emph{Scalar curvature distinguishes points with different local geometries. On smooth surfaces, such as the images shown on the left, points with positive/zero/negative curvature are locally like the surface of a sphere/flat plane/saddle point. On the right is illustrated how the curvature sign influences the local geometry around a point, in terms of triangles around the point (with different angle sums) and in terms of volume growth (where geodesic discs can grow faster/slower than in Euclidean space). The second row on the left shows some graphs which are naturally associated to spaces of constant curvature, see Section \ref{SS4.1: constant curvature graphs}.}}
    \label{fig: continuous curvature & constant curvature}
\end{figure}

The classical setting of curvature is differential geometry \cite{jost_riemannian_1995, berger_differential_1988}, where it characterizes how much and in which ways a smooth space differs from being flat around a point. \emph{Scalar curvature} for instance associates a single number to each point on a manifold and allows to distinguish the qualitatively different cases of positive curvature (like the surface of a sphere), zero curvature (like the flat plane) and negative curvature (like a saddle point), as illustrated in Figure \ref{fig: continuous curvature & constant curvature}. The scalar curvature is a measure of how much the volume of (small) $\epsilon$-balls around a point differ from the volume of $\epsilon$-balls in Euclidean space of the same dimension. \emph{Ricci curvature}, on the other hand, associates a tensor to each point on a manifold, and reflects the difference of volume growth between geodesics emanating from the point in two tangential directions compared to Euclidean growth. While the Ricci curvature includes `directional' information not present in the scalar curvature, the latter can be retrieved as the trace of the former. We refer the readers to \cite{jost_riemannian_1995, berger_differential_1988} for references on curvature and differential geometry.\\
While discrete spaces such as graphs lack the differential structure required for classical definitions of curvature, it is still widely acknowledged that these spaces exhibit relevant geometric features and consequently, there have been many attempts at formulating a theory of \emph{discrete curvature} that works in these generalized settings as well \cite{chung_logarithmic_1996, forman_bochners_2003, bakry_diffusions_2006, ollivier_survey_2010, saucan_simple_2020, Bobenko_2008_discrete}. The idea of discrete curvature is perhaps most easily understood in the case where the discrete spaces under consideration are closely related to some smooth space, for instance a graph which can be isometrically embedded into a manifold (see Section \ref{SS4.1: constant curvature graphs}) or random geometric graphs which are natural `coarse' representations of some smooth space (see Section \ref{SS4.3: resistance curvature for ERGs}). However, even when there is no clear underlying or latent smooth space involved, discrete notions of curvature can be formulated. Importantly, most contributions on the subject of discrete curvature not only attempt to generalize the definition of curvature, but also try to translate the multitude of curvature-related results over to this new setting. For more information on different approaches to
discrete curvature and their theory, we refer the reader to \cite{kamtue_combinatorial_2018, peyerimhoff_curvature_nodate}\cite{Najman_2017_modern, Bobenko_2008_discrete}. In Section \ref{SS4.2: relation to other discrete curvatures} we give a short introduction to combinatorial curvature, Ollivier-Ricci curvature and Forman-Ricci curvature, as these particular notions are related to our proposed curvature.
%%%%

%%%%
\subsection{Graphs, Laplacians and the effective resistance}
Here, we introduce the preliminaries on graphs and effective resistances required to define the resistance curvature; further details can be found in \cite{devriendt_effective_2020}. While we mainly work in the setting of \emph{weighted graphs}, it is important to note that a weighted graph is nothing more than a set together with some notion of similarity defined between (some) pairs of elements of this set. This type of data is of course very general and widely available, which broadens the applicability of the resistance curvature to more than just graph theory and network science. For instance, any dataset with a similarity kernel -- or a parametrized family of kernels to cover different scales, as is often the setting in topological data analysis \cite{ghrist_barcodes_2008, carlsson_topology_2009} -- will naturally have an associated weighted graph whose resistance curvature then says something about the geometry of the underlying dataset. This is illustrated by the analysis in Section \ref{SS4.3: resistance curvature for ERGs}, where we find that the resistance curvature (calculated via an associated graph) retrieves the flat curvature of a point-set in the Euclidean plane.
\\~\\
A weighted graph $G=(\mathcal{N},\mathcal{L},c)$ consists of a set of $n$ nodes $\mathcal{N}$ and a set of $m$ links $\mathcal{L}$ that connect pairs of nodes\footnote{Traditionally, nodes and links are often referred to as vertices and edges instead. We choose to work with ``nodes and links'' to describe graphs and reserve the terminology of ``vertices and edges'' to describe simplices which appear in the study of effective resistances, see for instance \cite{fiedler_matrices_2011, devriendt_effective_2020}. To avoid confusion, we also note that our notion of link is distinct from the ``link of a simplex'' in the context of simplicial complexes}, denoted by an unordered tuple $(i,j)\in\mathcal{L}$ or simply $i\sim j$. Furthermore, each of the links is assigned a positive weight through the function $c:\mathcal{L}\rightarrow \mathbb{R}_{>0}$ and we write $c_{ij}$ for the weight of a link $(i,j)$. We furthermore restrict our attention to \emph{simple, undirected graphs without self loops}. 
\\
One convenient way to describe and study weighted graphs is by their associated Laplacian matrix. For a given $n$-node graph, the \emph{Laplacian matrix} is an $n\times n$ matrix $Q$ whose rows and columns are identified with the node set $\mathcal{N}$, and with entries
$$
(Q)_{ij} := \begin{cases}
-c_{ij} \text{~if $i\sim j$}\\
k_i \text{~if $i=j$}\\
0\text{~otherwise},
\end{cases}
$$
where $k_i:= \sum_{j\sim i}c_{ij}$ is the \emph{weighted degree} of node $i$, i.e. the total weight of all links incident on node $i$. We will also be interested in the \emph{combinatorial degree} of a node $i$, which is defined as $d_i := \sum_{j\sim i}1$ and is independent of the link weights. Apart from being a convenient matrix representation of the data associated to a weighted graph, the Laplacian matrix is also a central object in (spectral) graph theory \cite{chung_spectral_1997, mohar_laplacian_1991, merris_laplacian_1994} with broad applications \cite{hoory_expander_2006, vishnoi_lxb_2012, spielman_graphs_2016, coifman_diffusion_2006, von_luxburg_tutorial_2007}. 
\\
The Laplacian matrix is positive semidefinite and its rank depends on the number of connected components $\beta(G)$ in the graph as \cite{chung_spectral_1997}
$$
\operatorname{null}(Q) = \beta(G) \text{~or~} \operatorname{rank}(Q) = n-\beta(G),
$$
where the kernel is spanned by vectors which are piecewise constant on the components of $G$; we write $G_i$ to denote the unique connected component that contains a given node $i$. In context of this article, the Laplacian is mainly relevant for its relation to the effective resistance. This resistance is defined based on the \emph{Moore-Penrose pseudoinverse} of the Laplacian $Q^\dagger$, which is determined by the equations $Q^\dagger Q = QQ^\dagger = \operatorname{proj}(\ker(Q)^\perp)$ and can be calculated, for instance, by inverting the nonzero eigenvalues of the Laplacian matrix.
\\
The \emph{effective resistance} between pairs of nodes is defined based on this pseudoinverse Laplacian as 
\begin{equation}\label{eq: definition effective resistance}
\omega_{ij} := (\mathbf{e}_i-\mathbf{e}_j)^TQ^\dagger(\mathbf{e}_i-\mathbf{e}_j),
\end{equation}
for any $i,j$ in the same connected component\footnote{Formally, one can take the effective resistance to be infinite between nodes in different components.}, where $\mathbf{e}_i$ is the $i^{\text{th}}$ unit vector. The \emph{resistance matrix} $\Omega$ is the $n\times n$ matix containing all pairwise effective resistances, as $(\Omega)_{ij}=\omega_{ij}$. The concept of effective resistances originates from the theory of electrical circuits, where it captures the resistance exerted by the whole network on a current flowing between any two nodes; more precisely, if a unit current flows through the graph from node $i$ to $j$, where it passes through the links with resistances $1/c$ and obeys Kirchhoff's laws, then the effective resistance is measured\footnote{The algebraic definition \eqref{eq: definition effective resistance} follows from this `measurement definition' because the linear relation between voltages $\mathbf{v}$ and (external) currents $\mathbf{x}$ imposed by Kirchhoff's laws can be expressed using the Laplacian matrix as $\mathbf{x}=Q\mathbf{v}$.} as the voltage difference between $i$ and $j$ (resistance$=$voltage/current) \cite{dorfler_electrical_2018, thomassen_resistances_1990}. The effective resistance has many mathematical properties and one of its most celebrated properties is that it is a metric between the nodes of a graph, see also \cite{devriendt_effective_2020}:
\begin{theorem}[Resistance is distance \cite{gvishiani_metric_1987, klein_resistance_1993}]\label{th: resistance is distance}
The effective resistance is a metric between the nodes of a graph.
\end{theorem}
\emph{Rayleigh's monotonicity principle} \cite{doyle_electric_1988, thomassen_resistances_1990} -- which says that adding new paths or shortening existing paths between two nodes can only decrease the effective resistance between these nodes, and vice versa for removing paths -- leads to the following intuitive interpretation of this distance measure:
Theorem \ref{th: resistance is distance} provides a good way to think intuitively about effective resistances: if the resistance distance $\omega_{ij}$ is large, this reflects that nodes $i$ and $j$ are not well connected in the network -- there are few and mainly long (low weight) paths connecting $i,j$ --, whereas a small $\omega_{ij}$ means that they are well connected -- there are many short paths between $i,j$. The effective resistance thus reflects a more integrated notion of distance compared to the shortest path distance, since it takes into account all the different paths between two nodes and how they are interconnected.

In this article, we will mainly be interested in the product of the link weight and the effective resistance between the end nodes of the link, $c_{ij}\omega_{ij}$ which we call the \emph{relative resistance}\footnote{This is not standard terminology, but is introduced here because it has distinctly different properties and interpretations than the effective resistance $\omega_{ij}$. The name relative resistance follows from $c_{ij}$ being a conductance (inverse resistance) in the context of electrical circuits, such that $c_{ij}\omega_{ij}$ is the quotient between a link's effective resistance and its `direct' resistance as a resistor.} of $(i,j)$. This relative resistance reflects the \emph{importance} of a link for connectivity of the graph. A link which is very redundant -- i.e. by connecting two nodes in group of tightly interconnected nodes -- will have a low relative resistance, whereas removing a link with high relative resistance would significantly reduce the connectivity between its endpoints. In \cite{spielman_spectral_2011}, it was shown that sampling links proportional to their relative resistance yields statistically representative sparse samples of a graph, which is related to the concept of statistical load/leverage in numerical linear algebra \cite{drineas_effective_2010}.

To further quantify the notion of `link importance' we introduce a well-known relation between relative resistances and random spanning trees -- since the relative resistance of a link only depends on the connected component containing the link we will briefly assume the underlying graph to be connected. A \emph{spanning tree} $T$ of $G$ is a connected subgraph on all the nodes which is a tree (i.e. has $n-1$ links) and the set of all spanning trees is denoted by $\mathcal{T}$. A \emph{random spanning tree} $\mathbf{T}$ is a random element of $\mathcal{T}$ with probability to equal any specific spanning tree $T$ given by
$$
\Pr[\mathbf{T}=T] = \frac{c(T)}{\sum_{T'\in\mathcal{T}}c(T')} \quad\quad\text{with}\quad\quad c(T):= \prod_{(i,j)\in\mathcal{L}(T)}c_{ij}
$$
In other words, a random spanning tree is sampled proportional to $c(T)$, the product of its link weights. This relates to the relative resistance as, see also \cite{lyons_determinantal_2003}:
\begin{theorem}[\cite{kirchhoff_ueber_1847, burton_local_1993}]\label{th: relative resistance and spanning trees}
The relative resistance of a link equals the probability that this link is contained in a random spanning tree, $c_{ij}\omega_{ij}=\Pr[(i,j)\in \mathbf{T}]$.
\end{theorem}
For graphs on multiple connected components, the same result holds in terms of random spanning forests which consist of a random spanning tree on each of the connected components independently. Theorem \ref{th: relative resistance and spanning trees} is an important tool to study the relative resistance, and we will make use of the following properties (see Appendix \ref{A: properties of relative resistance})
\begin{property}[relative resistance properties]\label{property: bounds for relative resistance and Foster's theorem}
The relative resistance of a link $(i,j)$ is bounded by $
0< c_{ij}\omega_{ij} \leq 1,$ with equality in the upper bound if and only if it is a cut link. The sum over all relative resistances equals $\sum_{(i,j)\in\mathcal{L}}c_{ij}\omega_{ij} = n-\beta$, which is known as Foster's Theorem.
\end{property}
%%%

%%%
\section{Introducing the resistance curvature}\label{S3: introducing resistance curvature}
We will define two types of discrete curvature: a \emph{scalar curvature} defined on the nodes of a graph as $p:\mathcal{N}\rightarrow\mathbb{R}$ and a \emph{Ricci curvature} defined on the links of a graph as $\kappa:\mathcal{L}\rightarrow\mathbb{R}$. We will refer to the nodal curvature simply as `resistance curvature' and will use the slightly longer `link resistance curvature' for the latter\footnote{While finishing the manuscript, we found that the name ``effective resistance curvature" was introduced in \cite{grippo_effective_2016} for a different notion of curvature derived from effective resistances.}. For further details and stronger versions of the results in this section, we refer to Appendices \ref{A: properties of node resistance curvature p} and \ref{A: properties of link resistance curvature kappa}.
\begin{definition}[resistance curvature]
The resistance curvature is defined as
\begin{equation}\label{eq: definition resistance curvature}
p_i := 1 - \frac{1}{2}\sum_{j\sim i}\omega_{ij}c_{ij} \textup{~for any node $i\in\mathcal{N}$}.
\end{equation}
\end{definition}
The resistance curvature is a node function, but instead of writing $p(i)$ we will mostly write $p_i$ for a node $i$ and use the $n\times 1$ vector $\mathbf{p}=(p_1,\dots,p_n)^T$ containing all resistance curvatures\footnote{To be fully consistent, we can assume a bijection between $\mathcal{N}$ and the integers $[n]=1,\dots,n$ such that we can index the nodes of $G$ both as $i\in\mathcal{N}$ or as $1\leq i\leq n$.}. The definition of $p$ appeared before in \cite{zhou_resistance_2016, bapat_resistance_2004}, and \cite[Proposition I.2]{devriendt_effective_2020} presents a number of alternative expressions for $p$ (see also Appendix \ref{A: alternative definitions}), but to our knowledge the connection to curvature has not yet been made. Initially, this connection is suggested by the following qualitative observation: if the neighbourhood of a node is tree-like (resp. clustered) with few (many) short cycles, then the local relative resistances will be large (small) and thus $p_i$ will be small (large), as expected for a notion of curvature \cite[Thm. 1]{Jost_Liu_2013_Ollivier}. Section \ref{SS: average distance characterization} introduces an alternative definition for the resistance curvature which further supports this intuition.
\\
While the resistance curvature is defined \emph{locally} as a sum over the relative resistances of the incident links, the relative resistance and thus $p_i$ depend on the structure of the whole graph in general. Thanks to efficient algorithms for Laplacian systems however, the resistance curvature can still be approximated efficiently, as discussed at the end of this section.
\\
Following Property \ref{property: bounds for relative resistance and Foster's theorem} for the relative resistance, we can find some first basic results on the node resistance curvature (see Appendix \ref{A: properties of node resistance curvature p}):
\begin{property}\label{property: properties for p_i}
The resistance curvature of a node with $d_i\geq 1$ is bounded by $1-d_i/2\leq p_i\leq 1/2$, with equality in the lower bound if and only if all links $j\sim i$ are cut links. The sum over all resistance curvatures equals $\sum_{i\in\mathcal{N}}p_i=\beta$.
\end{property}
In tree graphs, every link is a cut link and thus Property \ref{property: properties for p_i} implies that the non-leaf nodes in a tree have nonpositive curvature. For path graphs -- which are tree graphs with nodes of degree $1$ and $2$ -- we find that the two end nodes have resistance curvature $1/2$ while all other nodes have zero curvature. The bounds for $p$ in Property \ref{property: properties for p_i} thus produce some examples where negative and zero curvature happen generically, i.e. in tree and path graphs. To add an example of positive curvature, we can consider the cycle graph. By virtue of the rotational symmetry of the cycle graph, all nodes are indistinguishable and should thus have the same resistance curvature; by $\sum p_i=1$ this then implies that the resistance curvature is positive and equal to $p_i=1/n$ for all nodes in the $n-$cycle graph (see Appendix \ref{A: res curvature and transitivity} for more detail). These three examples are treated in further detail at the end of this section.
\\~\\
Starting from the (scalar) node resistance curvature we now introduce the (Ricci) link resistance curvature:
\begin{definition}[link resistance curvature] The link resistance curvature is defined as
\begin{equation}\label{eq: link resistance curvature definition}
\kappa_{ij} := \frac{2(p_i+p_j)}{\omega_{ij}} \textup{~for any link $(i,j)\in\mathcal{L}$}.
\end{equation}
\end{definition}
The link resistance curvature equals the sum of the (nodal) resistance curvature of its end nodes, divided by their effective resistance. Consequently, calculating the link resistance curvature is straightforward once $p$ is known, and when $p_i$ and $p_j$ have the same sign this immediately determines the sign of $\kappa_{ij}$ as well. Tree, path and cycle graphs can thus serve as a first example of graphs with links of negative/zero/positive link curvature.
\\
Following the bounds for $p$ we find the following bounds for $\kappa$ (see Appendix \ref{A: properties of link resistance curvature kappa} for a stronger version):
\begin{property}\label{property: bounds for kappa_ij}
The link resistance curvature is bounded by
$(4-d_i-d_j)/\omega_{ij}\leq \kappa_{ij} \leq 2/\omega_{ij}$, with equality in the lower bound if and only if all incident links to $i$ and $j$ are cut links.
\end{property}
Both $p$ and $\kappa$ can be approximated efficiently. While computing effective resistances based on definition \eqref{eq: definition effective resistance} requires a costly matrix inversion, the specific structure of Laplacian matrices can be exploited to approximate the effective resistance of any link in the graph with accuracy $(1-\epsilon)$ in time $O\left(m\log(c_{\max}/c_{\min})/\epsilon^2\right)$, i.e. basically in linear time in the number of links \cite{spielman_spectral_2011, vishnoi_lxb_2012}. As a consequence, $p$ and $\kappa$ can be approximated efficiently.
\\~\\
Our terminology for $p$ and $\kappa$ as scalar and Ricci curvature follows from the relation to existing notions of discrete scalar and Ricci curvature discussed in Section \ref{SS4.2: relation to other discrete curvatures}. However, in the continuous case these two types of curvature satisfy specific relations which are not reproduced in the context of $p$ and $\kappa$ -- for instance, scalar curvature is the trace of the Ricci curvature tensor while we have not found a similar sum formula that holds for resistance curvatures. This suggests that our choice of terminology might not be fully consistent and that perhaps $p$ or $\kappa$ is only an approximate or derived notion of discrete curvature.
\subsection{Examples} \label{SS3.1: examples for resistance curvature}
\begin{figure}[h!]
    \centering
    \includegraphics[width = 0.85\textwidth]{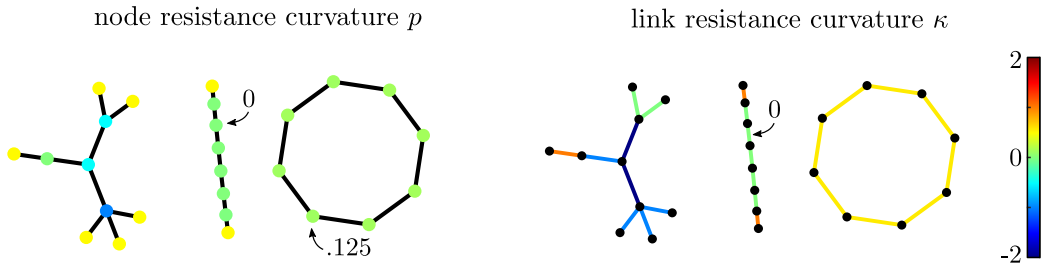}
    \caption{The node and link resistance curvatures in a tree, path and cycle graph (left to right).}
    \label{fig: path tree cycle}
\end{figure}

\begin{figure}[h!]
    \centering
    \includegraphics[width=\textwidth]{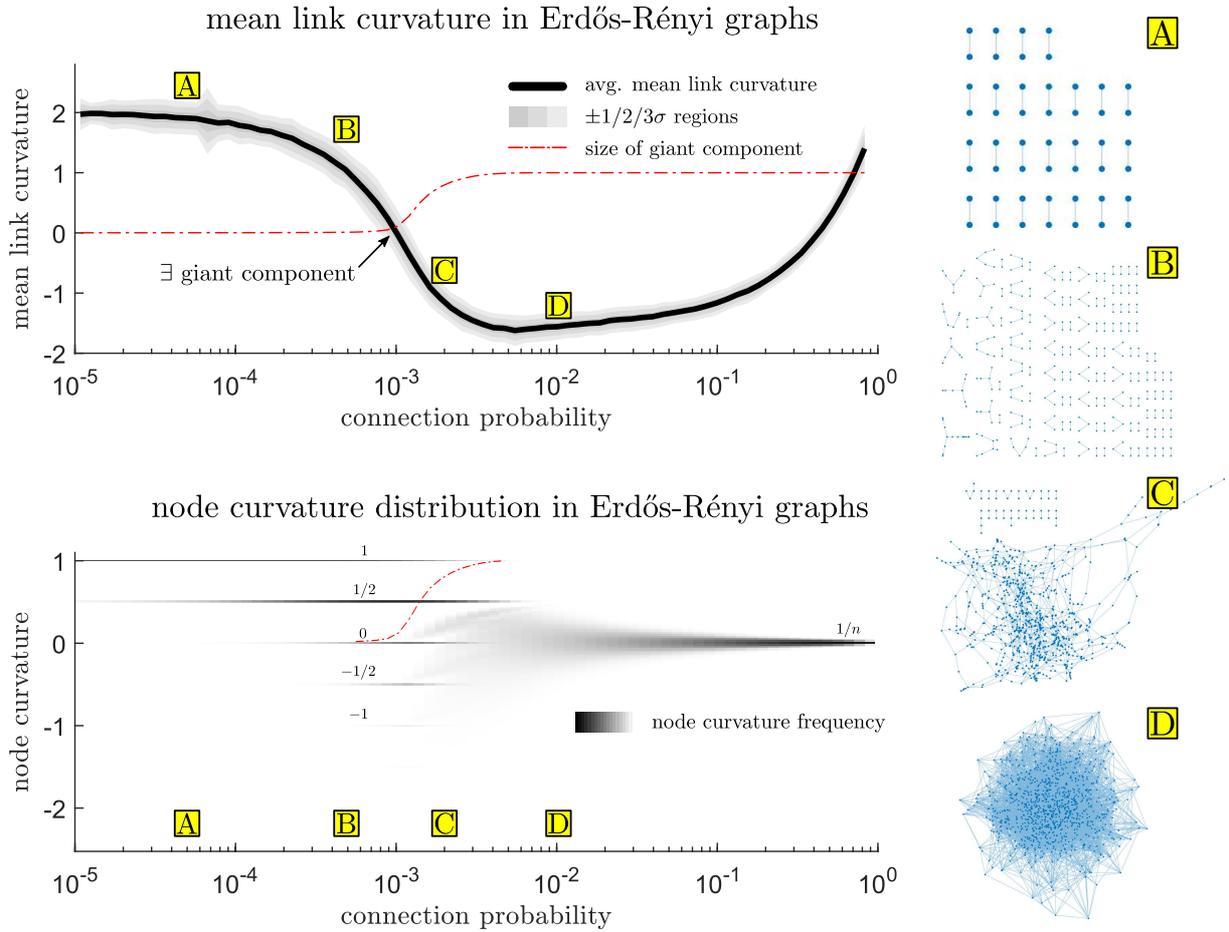}
    \caption{\emph{Evolution of the mean link resistance curvature and node curvature distributions in Erd\H{os}-R\'{e}nyi (ER) graphs for different connection probabilities (densities) $\rho$. We select $N_\rho = 10^4$ connection probabilities $\log_{10}(\rho)\in[-5,-0.05]$ uniformly at random, construct an ER random graph on $n=10^3$ nodes for each probability and calculate the mean link resistance curvature $\tfrac{1}{m}\sum_{j\sim i}\kappa_{ij}$ and node resistance curvatures $p$ in each graph. \textbf{\textup{The top left figure}} shows the sample average of the mean link resistance curvature, with graph samples binned together according to their connection probabilities ($N_{\text{bins}}=75$), and $\pm 1/2/3$ standard deviations around the average. In red, the size of the giant component ($=$ no. nodes in largest component/n) is plotted, again averaged according to the bins; we remark that the appearance of a giant component around $(n-1)\rho=1$ coincides with the mean link curvature crossing zero.  \textbf{The bottom left figure} shows the aggregate distribution of node curvatures for each of the $75$ connection probability bins. The shaded regions are proportional to the observed frequency of nodes with a given curvature $p$ in ER graphs in a given density bin $\rho$ -- white-black corresponds to $0$-$100\%$ observed frequency and to enhance the contrast (emphasize the lower frequencies), we plot the fourth power $f^4$ of the frequency $f\in[0,1]$.} \textbf{\textup{The right figure}} shows four ER graphs (A-D) with different connection probabilities to illustrate the evolution from a collection of many small tree-like components (A,B) for small $\rho$ to the emergence of a larger locally tree-like component (C) and finally, a dense connected graph (D) as $\rho$ further increases -- components consisting of single nodes are omitted for clarity}
    \label{fig: Erdos-Renyi}
\end{figure}

We discuss a few examples to support the introduced definitions; we describe some classes of graphs where positive, zero and negative curvature occurs generically and show some results for Erd\H{o}s-R\'{e}nyi random graphs. Figures \ref{fig: path tree cycle} and \ref{fig: Erdos-Renyi} illustrate the examples.
\\
\textbf{Example 1 (tree graphs)} As introduced, a tree graph is a connected graph with $m=n-1$ links. Since every link in a tree graph is a cut link, the bounds in Properties \ref{property: bounds for relative resistance and Foster's theorem}-\ref{property: bounds for kappa_ij} hold with equality, and we find that
$$
c_{ij}\omega_{ij}=1\quad\quad p_i=1-\frac{d_i}{2}\quad\quad \kappa_{ij} = c_{ij}(4-d_i-d_j).
$$
for all nodes $i$ and links $(i,j)$. In other words, except for the leaf nodes and the links connected to them, all nodes and links have nonpositive node/link resistance curvature.
\\
\textbf{Example 2 (path graphs)} A path graph is a tree with two leaf nodes forming the ends of the path, and all other nodes having degree two. Following the result for tree graphs we find that
$$
\begin{cases}
p_i = \frac{1}{2} \text{~at the end nodes}\\
p_i = 0 \text{~otherwise}
\end{cases} \quad\quad
\begin{cases}
\kappa_{ij} = c_{ij} \text{~if $i$ or $j$ is an end node}\\
\kappa_{ij} = 0 \text{~otherwise}
\end{cases}
$$
In other words, except for the end nodes and the links connected to them, all nodes and links have zero curvature. The case of a $2$-node path is slightly different as we get $\kappa_{ij}=2c_{ij}$.
\\
\textbf{Example 3 (cycle graphs)} As noted before, the symmetries of a graph can have implications for the resistance curvature. In Appendix \ref{A: res curvature and transitivity} we show that the relevant notion of symmetry is node (link) \emph{transitivity} which says that any node (link) can be mapped onto any other node (link) by some structure-preserving map, i.e. an automorphism. One example of a node and link transitive graph is the cycle graph, for which the node and link resistance curvatures equal (see Appendix \ref{A: res curvature and transitivity}):
$$
c\omega_{ij} = \frac{n-1}{n}\quad\quad p_i = \frac{1}{n}\quad\quad \kappa_{ij} = \frac{4c}{n-1}
$$
for all nodes $i$ and links $(i,j)$ and constant link weight $c$. In other words, all nodes and links are positively curved. Some more examples of transitive graphs will be discussed in Section \ref{SS4.1: constant curvature graphs}. Figure \ref{fig: path tree cycle} shows an example of tree, path and cycle graphs and their curvatures.
\\
\textbf{Example 4 (Erd\H{o}s-R\'{e}nyi random graph)} An \emph{Erd\H{o}s-R\'{e}nyi} (ER) graph is a random graph where every pair of nodes is connected with probability $\rho$. Figure \ref{fig: Erdos-Renyi} shows how the mean link curvature $\tfrac{1}{m}\sum_{i\sim j}\kappa_{ij}$ and resistance curvature distributions evolve in these random graphs as a function of the connection probability (density) $\rho$. An intuitive explanation for the observed evolution of the mean link curvature is as follows: ER graphs with very small $\rho$ typically consist of a collection of $2$-paths and disconnected nodes (graph $A$ in Figure \ref{fig: Erdos-Renyi}) with a mean link curvature of $+2$ since $\kappa_{ij}=2$ for a $2$-path. For increasing $\rho$, the ER graph will be a collection of trees or locally tree-like components (graphs $B,C$) resulting in a decreasing $\kappa$, and eventually the graph will become a dense connected graph (graph $D$) with increasing $\kappa$ until it reaches mean link curvature $+2$ for the complete graph at $\rho=1$. For the node resistance curvature we find a similar explanation: small and medium-density ER graphs consist mainly of disconnected nodes and $2$-paths with $p=1$ and $p=1/2$ respectively, or a collection of trees with $p_i=1-d_i/2$, which results in the `discrete distribution' of $p$ on half integers as observed in Figure \ref{fig: Erdos-Renyi}. For increasing density $\rho$, we observe that the locally tree-like graphs have a `continuous distribution' of possible node curvatures which narrows and converges to $p\approx 1/n$ when the density increases to $1$; this can be explained by the high uniformity of dense ER graphs in combination with Foster's Theorem for $\beta=1$. In Section \ref{SS4.3: resistance curvature for ERGs} we treat Euclidean random graphs as another example of random graphs.
%%%%%
\subsection{Average distance characterization}\label{SS: average distance characterization}
Definitions \eqref{eq: definition resistance curvature} and \eqref{eq: link resistance curvature definition} for the node and link resistance curvature are fairly simple to state and easy to manipulate, but they lack a geometric character that intuitively relates them to curvature. Here, we discuss an alternative definition of the resistance curvatures in terms of resistance distances in the neighbourhood around a node or link.
\\
A \emph{distribution} on a graph is a nonnegative function on the nodes $f:\mathcal{N}\rightarrow[0,\infty)$ with unit sum $\sum_{i\in\mathcal{N}} f(i)=1$. A distribution determines a random node $N\sim f$, which is a random element of the node set with probability of being equal to any particular node given by $\Pr[N=i]=f(i)$; we also say that $N$ is distributed according to $f$. The \emph{average distance between two random nodes $N$ and $M$}, distributed according to $f$ and $g$, can be defined as
\begin{equation}\label{eq: average distance}
\mathbb{E}(\omega_{NM}) := \sum_{i,j\in\mathcal{N}}f(i)g(j)\omega_{ij} \text{~with $N\sim f$ and $M\sim g$}.
\end{equation}
An important class of distributions follows from the diffusion (heat) equation on graphs\footnote{Alternatively, the normalized diffusion equation can be used $d\mathbf{f}(t)/dt=-Q\diag(k)^{-1}\mathbf{f}$; this lead to $p_i/k_i$ and $2(p_i/k_i+p_j/k_j)/\omega_{ij}$ in Proposition \ref{proposition: distance definitions for resistance curvature}, as also discussed in Section \ref{SS4.2: relation to other discrete curvatures}.}: $d\mathbf{f}(t)/dt=-Q\mathbf{f}$, where $Q$ is the graph Laplacian. If the initial state $\mathbf{f}(0)$ is a distribution then its evolution $\mathbf{f}(t)=\exp(-Qt)\mathbf{f}(0)$ remains a distribution; $\mathbf{f}(t)$ can also be interpreted as the distribution of a continuous-time random walker with initial distribution $\mathbf{f}(0)$. As we will discuss in more detail in Section \ref{SSS 2: OR curvature}, the diffusion of a distribution on a single node $f=e_i$ is concentrated on the neighbourhood of this node for small times; we will denote this by $\brho_{t,i}:=\exp(-Qt)\mathbf{e}_i$ which for small times satisfies $\brho_{t,i}\approx (I-Qt)\mathbf{e}_i$. In Appendix \ref{AA: proof of average distance characterization} we show that the node and link resistance curvature satisfy the following alternative definitions:
\begin{proposition}\label{proposition: distance definitions for resistance curvature}
The node and link resistance curvature are equal to
\begin{align}
&p_i = \lim_{t\rightarrow 0}\left(1 - \frac{1}{4t}\mathbb{E}(\omega_{N_tM_t})\right)\textup{~with $N_t,M_t\sim \brho_{t,i}$ independently}\label{eq: distance definition node resistance curvature}
\\
&\kappa_{ij} = \lim_{t\rightarrow 0}\frac{1}{t}\left(1 - \frac{\mathbb{E}\left(\omega_{N_tM_t}\right)}{\omega_{ij}}\right) \textup{~with $N_t\sim\brho_{t,i},M_t\sim\brho_{t,j}$}\label{eq: distance definition link resistance curvature}
\end{align}
\end{proposition}
In other words, the node resistance curvature is related to the \emph{average distance between nodes in its neighbourhood} and the link resistance curvature is related to \emph{the average distance between the neighbourhoods of the end nodes of the link}. Equivalently, this corresponds to the average distance between two independent random walkers (as defined by the diffusion equation) starting at a node, and to the average distance between two independent random walkers starting at the end nodes of a link. In Section \ref{SSS 2: OR curvature} we will show that expression \eqref{eq: distance definition link resistance curvature} is closely related to the definition of Ollivier-Ricci curvature. In \cite{devriendt_variance_2021} we found that $\var(\mathbf{f}):=\tfrac{1}{2}\mathbf{f}^T\Omega\mathbf{f}$ can be interpreted as the variance of a distribution on a graph; this allows to rewrite \eqref{eq: distance definition node resistance curvature} as
\begin{equation}\label{eq: curvature as variance}
p_i = \lim_{t\rightarrow 0} \left(1 - \frac{1}{2t}\var(\brho_{t,i})\right).
\end{equation}
The diffusion variance thus grows as $2t$ around nodes with zero node resistance curvature. This is equal to the variance growth for diffusion or, equivalently, the mean square displacement of Brownian motion on the real line (of zero curvature) \cite{Einstein_1905_brownian_motion}.
\\
Proposition \ref{proposition: distance definitions for resistance curvature} provides a (geo)metric interpretation of the resistance curvatures and in Appendix \ref{A: alternative definitions} we introduce some further alternative characterizations. In the remainder of this article we focus on definitions \eqref{eq: definition resistance curvature} and \eqref{eq: link resistance curvature definition} and their implications, but a further interpretation of Proposition \ref{proposition: distance definitions for resistance curvature} and the definitions in Appendix \ref{A: alternative definitions} and their consequences would be an interesting line of future research.
%%%%%
%%%%%

%%%%%
\section{Curvature properties of $p$ and $\kappa$}\label{S4: curvature properties}
In the previous section, we introduced the node and link resistance curvature based on effective (and relative) resistances. However, so far the only hint towards an interpretation of $p$ and $\kappa$ as notions of discrete curvature were some examples of graphs with negative/zero/positive curvature. In this section, we present three arguments that further evidence this interpretation. A brief summary of the arguments is given below and each argument is then presented in more detail in the following subsections. While these arguments do not constitute a `proof' that $p$ is a discrete scalar curvature and $\kappa$ a discrete Ricci curvature, we believe that they make a strong case and hope that this will provide a good starting point for further research. The arguments can be summarized as follows:
\begin{itemize}
    \item[(1)] The resistance curvatures are constant and of the correct sign for a number of graphs associated to \textbf{constant curvature spaces}. We show this for infinite regular lattices ($0$), infinite regular trees ($<0$) and platonic graphs as regular tilings of the sphere ($>0$).
    
    \item[(2)] The resistance curvatures are \textbf{related to established notions of discrete curvature} on graphs. We show that the node resistance curvature is related to combinatorial curvature and that the link resistance curvature satisfies `Forman curvature $\leq $ link resistance curvature $\leq$ Ollivier curvature'.
    
    \item[(3)] We show numerical and theoretical evidence that the node resistance curvature retrieves \textbf{zero curvature in Euclidean random graphs}.
\end{itemize}
\subsection{Constant curvature graphs}\label{SS4.1: constant curvature graphs}
\textbf{Positive curvature (platonic graphs)} Platonic graphs are the graph skeletons of the platonic solids (tetrahedron, cube, octahedron, dodecahedron, icosahedron), which are regular tilings of the $2$-sphere $\mathbb{S}^2$. By transitivity of these graphs, we know from Appendix \ref{A: res curvature and transitivity} that the resistance curvatures will be constant and positive, equal to $p_i=1/n$ and $\kappa_{ij}=2\rho$, where $\rho=m/{{n}\choose{2}}$ is the link density. This is in correspondence with the constant positive curvature of $\mathbb{S}^2$.
\\
\textbf{Zero curvature (infinite lattice)} The rectangular/triangular/hexagonal lattices are infinite graphs that correspond to (the graph skeleton of) the regular tilings of the Euclidean plane. The effective resistance in these lattices has been studied extensively and it was shown that the relative resistance of all links is equal to $2/d$, with $d$ the degree of the respective lattice \cite{doyle_electric_1988, thomassen_resistances_1990, flanders_infinite_1972}. Consequently, the lattice graphs have constant zero resistance curvatures $p_i=0$ and $\kappa_{ij}=0$, in correspondence with the zero curvature of $\mathbb{R}^2$.
\\
\textbf{Negative curvature (infinite regular tree)} The infinite regular tree or Bethe lattice \cite{mosseri_r_bethe_1982} corresponds to a regular tiling of the hyperbolic plane $\mathbb{H}^2$. Like all trees, the relative resistance of every link is equal to one which means that these graphs have negative resistance curvatures $p_i= 1-d/2$ and $\kappa_{ij}=2c(2-d)$ when degree $d>2$, in correspondence with the constant negative curvature of $\mathbb{H}^2$.
\\~\\
Some examples of these constant curvature graphs are shown in Figure \ref{fig: continuous curvature & constant curvature}. We remark that the correct curvature in regular tilings is not always reproduced correctly by the established notions of discrete curvature. As noted in \cite{kempton_large_2020} for instance, the hexagonal lattice is often assigned a negative curvature because it is locally tree-like.
%%%%%%%%%

%%%%%%%%%
\subsection{Relation to other discrete curvatures}\label{SS4.2: relation to other discrete curvatures}
As a second argument, we show that the node resistance curvature $p$ is related to combinatorial curvature (through random spanning trees) and that the link resistance curvature is related to Ollivier-Ricci and Forman-Ricci curvature (through tight two-sided bounds). 
\subsubsection{Node resistance curvature and combinatorial curvature}\label{SSS 1: combinatorial curvature}
Combinatorial curvature \cite{kamtue_combinatorial_2018, gromov_hyperbolic_1987} measures discrete curvature for graphs embedded in the plane (or other surfaces), i.e. with the nodes and links drawn in the plane such that the links only intersect at their endpoints. These drawings partition the plane into \emph{faces} which are connected regions of the plane bordered by links that form a cycle in the graph. The combinatorial curvature is then defined as \cite{kamtue_combinatorial_2018}
$$
p_i^{(co)} := 1-\frac{d_i}{2} + \sum_{\text{face~}f\ni i}\frac{1}{d_f},
$$
where $f\ni i$ are the faces that contain node $i$ in their boundary. For tree graphs, only a single unbounded face exists which does not contribute to the combinatorial curvature such that every node in a tree graph has combinatorial curvature $1-d_i/2$ equal to the node resistance curvature. Following the connection between relative resistances and random spanning trees, we moreover find that the resistance curvature of a node $i$ in a general graph equals the expected combinatorial curvature of $i$ in a random spanning tree; in other words (see Appendix \ref{A: correspondence with combinatorial curvature}):
\begin{equation}\label{eq: resistance curvature is expected combinatorial curvature} 
p_i = \mathbb{E}[p_i^{(co)}(\mathbf{T})] \text{~with random spanning tree $\mathbf{T}$}.
\end{equation}
\subsubsection{Resistance curvature and Ollivier-Ricci curvature}\label{SSS 2: OR curvature}
Yann Ollivier \cite{ollivier_survey_2010} introduced a notion of curvature for metric spaces with an associated Markov chain. This notion of discrete curvature can readily be applied to graphs where the shortest path or geodesic distance is a natural metric $d:\mathcal{N}\times\mathcal{N}\rightarrow \mathbb{R}$, and with lazy random walks taking the role of a Markov chain\footnote{The subscript `$t$' of the Markov chain is the \emph{laziness parameter} of the random walk and should not be confused with the `Markov time' of a continuous-time Markov chain which is often denoted by $t$ as well. For a given random walk with transition probabilities $T_{ij}=\Pr[j \text{~at step k+1} \vert i \text{~at step k}]$, the lazy random walk has transition probabilities $T'_{ij}=(1-t)I + t T_{ij}$, i.e. with probability $(1-t)$ to stay at the same node each step.} $\mu_t$. For every node $i$, the Markov chain determines a function $\mu_{t,i}:\mathcal{N}\rightarrow [0,1]$ on the nodes, which gives the probability to find the random walker at a node one step after being at $i$ (see later for some examples). This function can be thought of as a `ball' around $i$ in the graph. The metric $d$ can be lifted to provide a notion of distance between these balls by making use of the $1$-Wasserstein distance $W_1$:
$$
W_1(\mu_{t,i}, \mu_{t,j}) := \min_{P} \tr (PD),
$$
where the minimum is taken over all nonnegative $n\times n$ matrices $P$ with `marginals' $\sum_j(P)_{xj}=\mu_{t,i}(x)$ and $\sum_i(P)_{ix}=\mu_{t,j}(x)$ and where matrix $(D)_{ij}=d(i,j)$ is the \emph{distance matrix} of the graph. \emph{Ollivier-Ricci} (OR) curvature measures how much the \emph{direct distance} $d(i,j)$ between two nodes differs from the \emph{distance between balls} around these nodes $W_1(\mu_{t,i},\mu_{t,j})$ as a generalization of continuous Ricci curvature to the data $(\mathcal{N},d,\mu_t)$. Lin, Lu and Yau further modified Ollivier's definition as a limit for shrinking balls as \cite{lin_ricci_2011}
$$
\kappa^{(OR)}_{ij} := \lim_{t\rightarrow 0}\frac{1}{t}\left(1 - \frac{W_1(\mu_{t,i}, \mu_{t,j})}{d(i,j)}\right).
$$
For $t\rightarrow 0$, the distances $W_1$ and $d$ will converge and OR curvature thus measures how much they differ in the first order in $t$. If the distance between the balls is larger than the direct distance between the points then $\kappa^{(OR)}$ will be negative, and the other way around for positive curvature. 
\\~\\
One choice for Markov chain $\mu_t$ is the \emph{lazy random walk}
$$
\mu_{t,i}(j) = \begin{cases}
1-t k_i\text{~if $i=j$}\\
c_{ij}t \text{~if $j\sim i$}\\
0 \text{~else}
\end{cases}\text{~for $0\leq t\leq k_{\max}^{-1}$}
$$
which can be defined compactly as the vector $\bmu_{t,i} = (I-Qt)\mathbf{e}_i$, where $I$ is the identity matrix\footnote{For small $t$, the same random walk is obtained from the continuous-time Markov chain $\rho_{t,i}$ defined by the diffusion equation as in Section \ref{SS: average distance characterization}; this Markov chain was also used in \cite{gosztolai_unfolding_2021}.}. While the OR curvature has been studied with respect to this Markov chain and the geodesic distance, we have not found it being used in conjunction with the effective resistance distance. We find that OR curvature with respect to the data $(\mathcal{N},\omega,(I-Qt))$ relates to the link resistance curvature as (see Appendix \ref{A: correspondence with Ollivier curvature}):
\begin{proposition}\label{prop: bounds for OR curvature}
The Ollivier-Ricci curvature with respect to resistance distance and the lazy random walk is lower-bounded by the resistance curvature as $\kappa^{(OR)}_{ij} \geq \kappa_{ij}$, with equality if $(i,j)$ is a cut link.
\end{proposition}
In tree graphs, all links are cut links and the link resistance curvature will thus correspond to the OR curvature throughout the graph. Furthermore, we note that the resistance distance is equal to the geodesic distance in trees.
\\
An alternative Markov chain that is used more often is the normalized lazy random walk, defined by $\bmu_{t,i}=(I-Q\diag(\mathbf{k})^{-1}t)\mathbf{e}_i$ where $\diag(\mathbf{k})$ is the diagonal matrix with the weighted degrees on its diagonal \cite{lin_ricci_2011, munch_ollivier_2017}. Similar to the standard lazy random walk, we find the bound
\begin{equation}\label{eq: OR curvature and degree-normalized curvature}
\kappa^{(LLY)}_{ij} \geq \frac{2(p_i/k_i+p_j/k_j)}{\omega_{ij}}
\end{equation}
which suggests a degree-normalized version of the link resistance curvature. The superscript LLY refers to the authors Lin, Lu and Yau of \cite{lin_ricci_2011} where this variant of OR curvature was first considered.
\\~\\
\textit{Remark:} The definition of Ollivier-Ricci curvature, and in particular the Lin-Lu-Yau limit version \cite{lin_ricci_2011}, is closely related to definition \eqref{eq: distance definition link resistance curvature} as the leading term of the ratio of the distance between shrinking balls around the end nodes of a link and the direct distance between these nodes; especially since we may replace $\brho_{t,i}$ by $\bmu_{t,i}$ in \eqref{eq: distance definition link resistance curvature}. However, while the change from $\kappa_{ij}^{(OR)}$ to $\kappa_{ij}$ comes down to `just' changing shortest-path distance to resistance distance and Wasserstein distance $W_1$ to average distance, the resulting measure $\kappa_{ij}$ does seem to have a unique set of properties that cannot simply be explained by its relation to OR curvature. For instance, the alternative definitions in Appendix \ref{A: alternative definitions} and the associated discrete Ricci flow in Section \ref{S5: Discrete Ricci Flow} are not shared by the Ollivier-Ricci curvature and the example of hexagonal lattices, misclassified by OR curvature as negatively curved but assigned a zero curvature by the link resistance curvature illustrates that the two definitions disagree in important ways.
\subsubsection{Resistance curvature and Forman-Ricci curvature}\label{SSS 3: FR curvature}
In \cite{forman_bochners_2003}, Robin Forman introduced a new notion of curvature for CW complexes, a class of discrete/combinatorial spaces that includes graphs. This so-called Forman-Ricci (FR) curvature is defined by generalizing the definition of the classical Ricci curvature in terms of the Bochner Laplacian of a manifold to a definition for discrete spaces based on an analogous Bochner Laplacian on these spaces, see also \cite{jost_characterizations_2021}. The FR curvature is expressed in terms of \emph{local} combinatorial data around a considered point, and Sreejith et al. \cite{sreejith_forman_2016} translated Forman's general definition to graphs as
\begin{equation}\label{eq: definition FR curvature}
\kappa^{(FR)}_{ij} = 2w_i\left(1-\frac{1}{2}\sum_{k\sim i}\sqrt{\frac{w_{ij}}{w_{ik}}}\right) + 2w_j\left(1-\frac{1}{2}\sum_{k\sim j}\sqrt{\frac{w_{ij}}{w_{jk}}}\right)
\end{equation}
where $w_*$ are some nonzero weights associated to the nodes and links of the graph. Expression \eqref{eq: definition FR curvature} clearly resembles the definition of $\kappa$ in terms of the node resistance curvature and indeed, choosing unit weights yields the following relation (see Appendix \ref{A: correspondence with Forman curvature}):
\begin{proposition}\label{prop: upper bound for FR curvature}
The Forman-Ricci curvature with respect to unit weights is upper-bounded by the link resistance curvature as $\kappa^{(FR)}_{ij}/\omega_{ij} \leq \kappa_{ij}$,
with equality if and only if $(i,j)$ is a cut link.
\end{proposition}
\textit{Remark:} If we choose the inverse degree as node weights $w_i=1/k_i$ and consider graphs with $c=1$, then we find the bound $\kappa_{ij}^{(FR)}/\omega_{ij}\leq 2(p_i/k_i + p_j/k_j)/\omega_{ij}$ which again features the degree-normalized version of the link resistance curvature as in \eqref{eq: OR curvature and degree-normalized curvature}.
\\
\textit{Remark:} We note one further related notion of discrete curvature: in \cite{steinerberger_curvature_2022}, it is shown that for the shortest-path distance matrix $D$ the equation $D\mathbf{w}=\mathbf{u}$ determines a function $w$ on the nodes with many properties of a discrete curvature. This relates to resistance curvature by the expression $\Omega\mathbf{p}\propto \mathbf{u}$ and $\mathbf{u}^T\mathbf{p}=1$ introduced in Appendix \ref{A: alternative definitions}. Despite this close relation, the particular properties of the resistance matrix (for instance, its invertibility \cite{devriendt_effective_2020}) result in a solution $\mathbf{p}$ with different properties; for instance, most of the alternative definitions in Appendix \ref{A: alternative definitions} do not translate to $\mathbf{w}$.
\\~\\
\textbf{Summary} The relation between the link resistance curvature and Forman and Ollivier-Ricci curvature can be summarized by the two-sided bound
$$
\kappa_{ij}^{(OR)} \geq \kappa_{ij} \geq \kappa_{ij}^{(FR)}/\omega_{ij} \text{~with equality for cut links,}
$$
where we stress that these are OR and FR curvatures with respect to particular data: $\omega$ as metric for OR and $w=1$ as weight for FR. Apart from being an argument for the interpretation of the link resistance curvature as a discrete curvature, this result is relevant in the context of other works that investigate the relation between both curvatures \cite{jost_characterizations_2021, tee_enhanced_2021, samal_comparative_2018}.
\subsection{Resistance curvature for Euclidean random graphs}\label{SS4.3: resistance curvature for ERGs}
\begin{figure}[h!]
    \centering
    \includegraphics[width=0.85\textwidth]{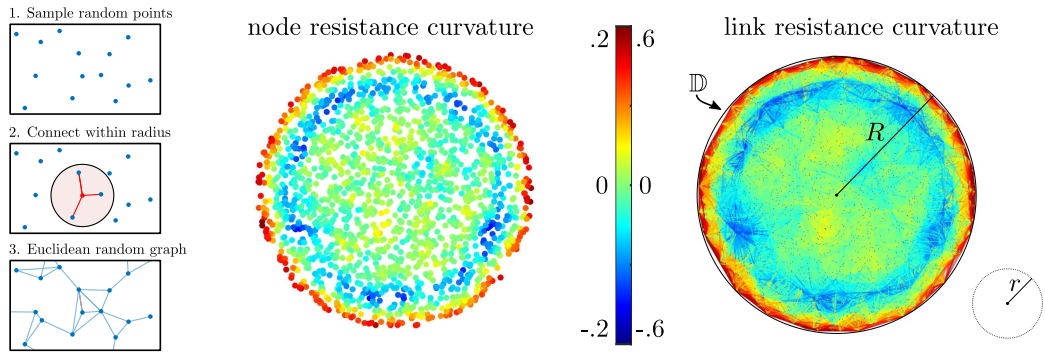}
    \caption{\emph{ \textbf{\textup{The left panels}} illustrate the construction of a Euclidean random graph (ERG): points are randomly sampled from a domain $\mathbb{D}\subseteq\mathbb{R}^2$ and connected if they lie within a certain connection radius $r$.  \textbf{\textup{The right figure}} shows the resistance curvature calculated in an ERG on a disc with radius $R/r=4$ and expected number of nodes $N=2.10^{3}$, with red/green/blue indicating positive/zero/negative values. Visually, both the node and link resistance curvatures are close to zero near the center of the disc but become negative and then positive going closer to the boundary. More detailed experimental results are shown in Figure \ref{fig: result for ERGs} and the origin of the boundary effect is explained in the main text and Appendix \ref{A: res curvature in ERGs}}.}
    \label{fig: constructing ERG and resistance curvature on disc}
\end{figure}

\begin{figure}[h!]
    \centering
    \includegraphics[width=\textwidth]{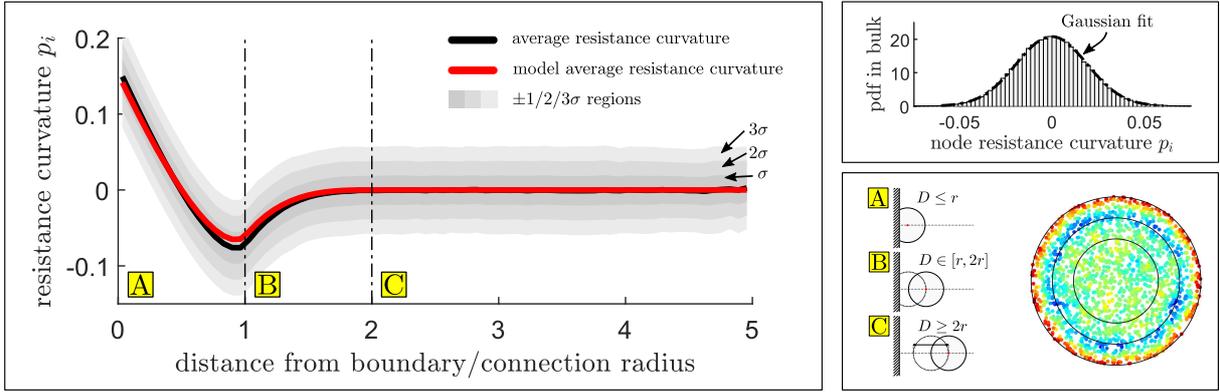}
    \caption{\emph{\textup{\textbf{Left panel:}} We sample $100$ graphs from an ERG on a disc with $R/r=5, N=5.10^3$ and calculate the node resistance curvature $p_i$ for all nodes. The plot shows the mean and variation around the mean of these curvatures (vertical axis) as a function of the distance from the boundary divided by the connection radius $D_i/r$ (horizontal axis); the mean and variation are calculated by grouping the $100\times N$ points according to their boundary distance in 75 equal-distance bins. The red line shows the heuristic $\mathbb{E}(\hat{p}(D))$ (expression \eqref{eq: boundary function formula}) which is in close correspondence to the black line, the experimental approximation for $\mathbb{E}(p(D))$. The boundary (A), near-boundary (B) and bulk regime (C) are delineated by vertical black lines and the \textup{\textbf{bottom right}} panel illustrates how these different regimes correspond to different local geometries of the nodes (see also Appendix \ref{A: res curvature in ERGs}). The \textup{\textbf{top right}} panel shows the distribution of $p_i$ for nodes in the bulk regime and the maximum likelihood Gaussian fit of this distribution; the good correspondence of this fit suggests that $p_i$ could be a zero-mean Gaussian random variable in the bulk.}}
    \label{fig: result for ERGs}
\end{figure}

A natural question for discrete curvatures is whether they converge to continuous curvature on discrete structures that represent finer and finer discretizations of some continuous space. For instance, van der Hoorn et al. \cite{van_der_hoorn_ollivier-ricci_2021} recently proved convergence of Ollivier-Ricci curvature to Ricci curvature in the continuum limit of random geometric graphs sampled from Riemannian manifolds. Here, as our third argument, we consider Euclidean random graphs and find numerical and theoretical evidence that the resistance curvature of these graphs retrieves the underlying zero curvature of the Euclidean plane.
\\
A Euclidean random graph (ERG) is a random graph constructed on some domain $\mathbb{D}\subseteq\mathbb{R}^2$ in the plane from which points are sampled uniformly at random by a Poisson point process with a homogeneous rate $\lambda$; the expected number of points is equal to $N=\lambda\text{Area}(\mathbb{D})$ so we may equivalently fix a desired $N$. These points are then taken as the nodes of a graph and pairs of nodes are linked if they lie within a certain connection radius (distance) $r$ from each other. All together, a Euclidean random graph model or ensemble can thus be parametrized by $(\mathbb{D},N,r)$. See \cite{penrose_random_2003} for more information on ERGs and their properties and Figure \ref{fig: constructing ERG and resistance curvature on disc} for an illustration of their construction and some examples.
\\
To simplify the setup for further analysis, we consider ERGs on the disc $\mathbb{D}=\lbrace x\in\mathbb{R}^2:\Vert x\Vert\leq R\rbrace$ of radius $R$, such that nodes at the same radial distance $D_i$ from the boundary are statistically equivalent, and the effective parameters reduce to $(R/r,N)$. 
\\
As a first numerical analysis, Figure \ref{fig: constructing ERG and resistance curvature on disc} shows the resistance curvatures in a Euclidean random graph on the disc with $R/r=4$ and $N=2.10^3$ nodes. Near the center of the disc, which we call the \emph{bulk} of the graph, we clearly observe that the curvatures are close to zero, while moving towards the boundary the resistance curvature becomes first negative and then positive at the boundary. Figure \ref{fig: result for ERGs} shows this result in more detail: the mean and standard deviations of $p_i$ are plotted with respect to $D_i/r$ with data aggregated over $100$ samples of an ERG with $R/r=5$ and $N=5.10^3$. The main observation in Figure \ref{fig: result for ERGs} is again that $p_i$ is close to zero in the bulk -- moreover, it appears to be a zero-mean random variable -- and becomes negative and then positive when moving closer to the boundary.
\\~\\
We now develop a model for the resistance curvature in Euclidean random graphs that partly explains these observations. While a full understanding of the relative resistance in ERGs is lacking, in particular concerning dependencies between the relative resistance of different links, we find that the numerical observations in Figures \ref{fig: constructing ERG and resistance curvature on disc} and \ref{fig: result for ERGs} can be explained in great detail using a simple heuristic for the effective resistance. In  \cite{luxburg_getting_2010, luxburg_hitting_2014}, von Luxburg et al. showed that in ERGs with increasing number of nodes, the effective resistance between any pair of nodes $i,j$ converges to the sum of their inverse degrees $d_i^{-1}+d_j^{-1}$. While the degree of a node depends on the specific random graph realization, we know by properties of Poisson sampling that the \emph{expected degree} of a given node $i$ is determined by the area of overlap between the domain and an $r$-radius disc centered at $i$, in other words $\mathbb{E}(d_i)=\lambda S(i)$ where $S(i):=\text{Area}(\lbrace x\in\mathbb{D} : \Vert x-i\Vert\leq r\rbrace)$, with expectation over the random graph ensemble. Combining these two results, we propose the following heuristic $\hat{\omega}$ for the effective resistance in ERGs:
\begin{equation}\label{eq: heuristic for resistance}
\hat{\omega}_{ij} := \frac{1}{\lambda S(i)} + \frac{1}{\lambda S(j)} \text{~for all links $i\sim j$.}
\end{equation}
Importantly, this reduces the effective resistance of a link to a purely \emph{geometric and local} quantity that is only determined by how the $r$-radius region around the link overlaps with the domain $\mathbb{D}$. In Appendix \ref{A: res curvature in ERGs}, we show that this heuristic combined with the assumption that $r\ll R$ allows to calculate the expected resistance curvature $\hat{p}$ as a function of the distance to the boundary $D$ as
\begin{equation}\label{eq: boundary function formula}
\mathbb{E}(\hat{p}(D)) = \frac{1}{2}\left(1 - \int_{t=-r}^{\min(D,r)}\frac{\vert dA(t)\vert}{A(\max(-r,t-D))}\right)\quad\quad
\raisebox{-0.45\height}{
\includegraphics[width=0.13\textwidth]{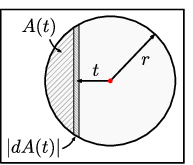}}
\end{equation}
where $A(t)$ is the area of a circular segment at height $t$ and $\vert dA(t)\vert$ the change of segment size -- as illustrated in the figure on the right -- and with expectation taken over the ERG ensemble. While we were not able to find a general closed-form expression for the integral in \eqref{eq: boundary function formula}, the expression can be evaluated numerically and is shown in Figure \ref{fig: result for ERGs} by the red line. Importantly, Figure \ref{fig: result for ERGs} shows that formula \eqref{eq: boundary function formula} matches very well with the experimentally observed resistance curvatures (the black line). This agreement suggests that our heuristic captures the main mechanisms behind the convergence of resistance curvature in Euclidean random graphs.
\\~\\
Due to the appearance of \emph{min} and \emph{max} operations, expression \eqref{eq: boundary function formula} is a piecewise function of the boundary distance $D$. In Appendix \ref{A: res curvature in ERGs}, we show that there are three possible regimes for $D$ as a result of the different local geometries around a node, i.e. how much of the connection discs around the node and its neighbours overlap with the domain (see also Figure \ref{fig: result for ERGs}). We find the following cases: (A) the \emph{boundary regime} $D\leq r$ where nodes as well as some of their neighbours are influenced by the boundary, (B) the \emph{near-boundary regime} $r\leq D\leq 2r$ where nodes are not influenced directly by the boundary (i.e. no overlap between the connection disc and the boundary), but some of its neighbours are, and (C) the \emph{bulk regime} $D\geq 2r$ where nodes are at least two connection radii away from the boundary such that the nodes nor any of their neighbours are influenced by the boundary. Most importantly, in the bulk regime expression \eqref{eq: boundary function formula} simplifies to $\mathbb{E}(\hat{p}(D))=0$ which means that we have {zero node resistance curvature $\hat{p}$ in expectation in the bulk of ERGs}. In the limit\footnote{While the ratio $r/R$ needs to be small to reduce the boundary effect, the expected number of nodes $N$ needs to be large enough with respect to $r$ such that the Poisson distribution (for the node degrees) concentrates around its mean and to be consistent with the conditions of Von Luxburg et al. for our heurstic \eqref{eq: heuristic for p for one sample}.} of $r/R\rightarrow 0$, this bulk regime will take up all but a vanishing fraction of the domain and almost all nodes will thus have zero expected resistance curvature.
\\
We remark that the derivation for expression \eqref{eq: boundary function formula} in Appendix \ref{A: res curvature in ERGs} is independent of the specific shape of the domain $\mathbb{D}$, i.e. it need not be a disc, and that for ERGs in higher dimensions $A(t)$ would simply be replaced by the volume of a higher-dimensional spherical cap.
\\~\\
Both the numerical results in Figure \ref{fig: constructing ERG and resistance curvature on disc} and the heuristic model suggest that the link resistance curvature also converges to zero (in the mean) in the bulk of Euclidean random graphs. However, the statistical dependencies between the node resistance curvatures at the ends of a link might obstruct convergence for graphs sampled from more general manifolds; a more detailed analysis of $p$ and $\kappa$ for random geometric graphs is an interesting line of future research.
%%%%%%

%%%%%%
\section{Discrete Ricci flow}\label{S5: Discrete Ricci Flow}
Ricci flow is an important concept associated to curvature in differential geometry \cite{chow_ricci_2004}. Intuitively, the Ricci flow describes a curvature-dependent evolution of the metric of a Riemannian manifold, where the metric decreases in directions of positive curvature while increasing in directions of negative curvature. The overall effect of this flow is that the curvature is `smoothed out' and for this reason, the Ricci flow is also thought of as a form of nonlinear heat equation for the metric. Due to its importance in the continuous setting, there have been several proposals to define a Ricci flow in the case of discrete curvatures, with various applications: Jin et al. \cite{jin_discrete_2008} define a Ricci flow for meshes in computer graphics, Ni et al. \cite{ni_community_2019} use a Ricci flow inspired evolution of networks for community detection and network alignment and Weber et al. considered Ricci flows related to the Forman Ricci curvature \cite{weber_forman-ricci_2016}. Recently, Bai et al. studied a Ricci flow based on the Lin-Yu-Yau curvature \cite{bai_ollivier_2021} and Cushing et al. studied a Ricci flow based on Bakry-\'{E}mery curvature \cite{Cushing_2022_flow}. Here, we show that there is a natural Ricci flow associated to the resistance curvature with promising features.
%%%%%%

%%%%%%
\subsection{Resistance Ricci flow}
Ollivier proposed in \cite[Problem N]{ollivier_survey_2010} to look at the differential equation $\tfrac{d}{dt}d(i,j)=-\kappa^{(OR)}_{ij}d(i,j)$ as a discrete version of the Ricci flow. If we introduce the link resistance curvature $\kappa_{ij}$ and resistance distance $d(i,j)=\omega_{ij}$ in this expression, we find
\begin{equation}\label{eq: definition resistance ricci flow}
\frac{d\omega_{ij}}{dt} = -2(p_i+p_j) \text{~for all $i\neq j$ in the same component}
\end{equation}
which we call the \emph{resistance Ricci flow}; since $\omega$ is metric, we implicitly have $d\omega_{ii}/dt=0$. Ideally, equation \eqref{eq: definition resistance ricci flow} would describe the evolution of a resistance matrix $\Omega(t)$ for all $t>0$ and starting from any initial resistance matrix $\Omega_0:=\Omega(0)$. Unfortunately however, $\Omega(t)$ is in general not guaranteed to be a resistance matrix as the flow can result in graphs with negative and diverging link weights in finite time as we show below. Nonetheless, we find that the resistance Ricci flow satisfies a number of interesting properties. As the flow \eqref{eq: definition resistance ricci flow} affects each connected component independently, we will assume $G$ to be connected in the rest of this section.
\\~\\
If $S$ is a state space (e.g. a vector space) and $P:S\rightarrow \mathbb{R}$ a potential defined on this state space, then the evolution $ds/dt=-\nabla P$ for a flow $s(t)\in S$, assuming that it is well-defined, is called a \emph{gradient flow}. This name reflects that the resulting flow follows the direction of steepest descent of the potential $P$, as given by the negative gradient. As a first result, we find that the resistance Ricci flow is a gradient flow (see Appendix \ref{A: res Ricci flow}):
\begin{proposition}\label{proposition: gradient expression for resistance Ricci flow}
The resistance Ricci flow \eqref{eq: definition resistance ricci flow} is a gradient flow of the potential $\tr(\tfrac{1}{2}\Omega Q\Omega)$ defined on symmetric zero-diagonal matrices $\Omega$, as
\begin{equation}\label{eq: gradient dynamics}
\frac{d\Omega}{dt} = -\nabla_{\Omega}\tr\left(\tfrac{1}{2}\Omega Q\Omega\right),
\end{equation}
\end{proposition}
\textit{Remark:} The gradient $\nabla_{\Omega}$ is defined on the state space of symmetric zero-diagonal matrices, such that the flow $\Omega(t)$ remains symmetric and with zero-diagonal; see Appendix \ref{A: res Ricci flow} for further details.
\\
The gradient form \eqref{eq: gradient dynamics} relates the resistance Ricci flow to the classical diffusion equation for a density $\mathbf{x}$ in $\mathbb{R}^n$, which is a gradient flow for the potential $\tr(\tfrac{1}{2}\mathbf{x}^TQ\mathbf{x})$. Expression \eqref{eq: gradient dynamics} could have practical implications for the resistance Ricci flow: since the potential is a decreasing function of time for gradient dynamics as $d\tr(\tfrac{1}{2}\Omega Q\Omega)/dt<0$, this could be a starting point to study convergence of the resistance Ricci flow.
\\
In Appendix \ref{A: res Ricci flow}, we furthermore show that the potential can be rewritten as $\tr(\tfrac{1}{2}\Omega Q \Omega)=2n\sigma^2$ where $\sigma^2=\frac{1}{2}p^T\Omega p$ is a graph invariant that appears in several contexts in the study of effective resistances, see for instance \cite{devriendt_effective_2020} and \cite[Appendix I]{devriendt_variance_2021}. The gradient expression \eqref{eq: gradient dynamics} further adds to the theory of this invariant. As noted in \cite{devriendt_effective_2020}, the number $1/(2\sigma^2)$ also relates to the so-called \emph{magnitude} (introduced in \cite{leinster_magnitude_2013}) of the metric space $(\mathcal{N},\omega)$; this suggests to consider gradient flows of metric spaces with the magnitude as a potential function.
%%%%

%%%
\subsection{Flow of Laplacians}
While expression \eqref{eq: definition resistance ricci flow} does not provide much insight into how the resistance Ricci flow affects the graph structure, we can also express the flow in terms of Laplacian matrices (see Appendix \ref{A: res Ricci flow}):
\begin{proposition}\label{propos: Laplacian form for Ricci flow}
The resistance Ricci flow \eqref{eq: definition resistance ricci flow} is equivalent to the following flow of Laplacian matrices:
\begin{equation}\label{eq: dQ=2Q^2}
\frac{dQ}{dt} = 2Q\diag(\mathbf{p})Q \textup{,~for some initial Laplacian $Q(0)=Q_0$}.
\end{equation}
\end{proposition}
Proposition \ref{propos: Laplacian form for Ricci flow} is remarkable since it shows that the resistance Ricci flow -- which was defined based on the proposal of Ollivier -- corresponds to a simple flow of Laplacian matrices. By considering the off-diagonal entries of expression \eqref{eq: dQ=2Q^2}, we can write the evolution of individual link weights as:
\begin{equation}\label{eq: resistance Ricci flow for link weights}
\frac{dc_{ij}}{dt} = -2\sum_{\substack{k\sim i,j\\k\neq i,j}}c_{ik}c_{kj}p_k + 2c_{ij}\left(k_ip_i + k_jp_j\right) \text{~for all $i\neq j$},
\end{equation}
where the equation for the diagonal ($dk_i/dt$) follows from the conserved zero row sum. The change of link weights in the resistance Ricci flow is thus caused by two types of processes: there is one term for each shared neighbour $k$ of $i$ and $j$ which can increase or decrease the link weight based on its resistance curvature $p_k$, and a second term based on the degrees and curvatures of $i$ and $j$. Other works such as \cite{aoki_self-organization_2015} have arrived at similar dynamical rules when modeling dynamic network structures, with a balance between positive link-reinforcing terms and competing link-decaying terms. As a possible variation on the resistance Ricci flow, one might consider the dynamics $dQ/dt=2Q\diag(\mathbf{q})Q$ where $q$ is some function on the nodes; for instance, if $q$ depends on a process taking place on the graph such as in \cite{aoki_self-organization_2015}, this could model the codependent evolution of structure and processes. We now consider two examples of the resistance Ricci flow.
\\
\textbf{Example 1 (node transitive graphs)} The node resistance curvature in node transitive graphs is constant and equal to $p_i=1/n$. Furthermore, from \eqref{eq: definition resistance ricci flow} we find that $d\omega_{ij}/dt=-2/n$ is equal for all pairs of nodes such that transitivity is conserved by the flow and thus $p_i(t)=1/n$ as long as this is well-defined. As a consequence, the resistance Ricci flow for node transitive graphs simplifies to $dQ/dt = 2Q^2/n$. Diagonalizing this equation and solving the differential equation $dx/dt=2x^2\Rightarrow x(t)=x_0/(1-2tx_0)$ for each eigenvalue, we find the solution
$$
Q(t) = \left[I-\tfrac{2}{n}tQ_0\right]^\dagger Q_0 \text{~for $t<\frac{n}{2\mu_{\max}(Q_0)}$}.
$$
This flow \emph{diverges in finite time} for $t\rightarrow n/(2\mu_{\max})$. The resistance Ricci flow for node transitive graphs is formally very similar to a matrix flow studied in \cite{marvel_continuous-time_2011} as a model for structural balance in social networks, which considers the differential equation $dX/dt = X^2$ for a matrix $X$. The entries $(X)_{ij}$ of this matrix represent pairwise affinities or rivalries, depending on the sign, and finite-time divergence was interpreted in the context of the model as the onset of a structurally balanced configuration of social relations.
\\
\textbf{Example 2 (path graph)} All nodes in a path graph have zero node resistance curvature, except for the end nodes which have $p=1/2$. From equation \eqref{eq: resistance Ricci flow for link weights} for the evolution of link weights, we then find that the only change occurs for the end links, whose link weight evolves according to
$$
\frac{dc}{dt} = c^2 \Rightarrow c(t) = \frac{c_0}{1-c_0t} \text{~for $t<\frac{1}{c_0}$},
$$
and diverges to $+\infty$ for $t\rightarrow 1/c_0$. While this result may seem pathological, we can interpret the `infinite affinity' (which corresponds to $\omega\rightarrow 0$) as merging the end node and its neighbour into a single node -- this is a common operation in the theory of effective resistances. This analysis shows that the resistance Ricci flow of a path graph evolves by merging the end nodes of the path with their neighbours at $t=1/c_0$, resulting in a path of decreasing length until just a single node remains.
\\~\\
Finally, we remark that the `degree-normalized' link resistance curvature $2(p_i/k_i+p_j/k_j)/\omega_{ij}$, which naturally appears in correspondence with Ollivier and Forman Ricci curvature (see \eqref{eq: OR curvature and degree-normalized curvature}), also has an associated Ricci flow $d\omega_{ij}/dt=-2(p_i/k_i+p_j/k_j) \equiv dQ/dt=2Q\diag(\mathbf{p}/\mathbf{k})Q$. Initial results suggest that this Ricci flow has additional interesting properties.
%%%%%%

%%%%%%
\section{Conclusion}\label{S6: conclusion}
This article proposes a new approach to discrete curvature based on the effective/relative resistance, and we introduce the \emph{node resistance curvature} $p:\mathcal{N}\rightarrow \mathbb{R}$ and \emph{link resistance curvature} $\kappa:\mathcal{L}\rightarrow \mathbb{R}$ in definitions \eqref{eq: definition resistance curvature} and \eqref{eq: link resistance curvature definition} respectively. We provide a background on the effective and relative resistance, derive some basic theoretical results on $p$ and $\kappa$ and present a number of arguments that support the interpretation of $p$ and $\kappa$ as a curvature: they retrieve the correct (or expected) curvature for the discretization of some constant curvature spaces (\S \ref{SS4.1: constant curvature graphs}), the resistance curvatures are related to established notions of discrete curvature, combinatorial, Ollivier-Ricci and Forman-Ricci curvature (\S\ref{SS4.2: relation to other discrete curvatures}) and we provide numerical and theoretical evidence that $p$ retrieves zero curvature in expectation for Euclidean random graphs (\S\ref{SS4.3: resistance curvature for ERGs}). As an application, we show that the resistance curvatures have a natural associated Ricci flow and we discuss some of its properties.
\\
While a full understanding of the relation to continuous curvature is still lacking, the resistance curvatures complement the existing study of discrete curvatures and, by their relation to the effective resistance with a rich theory, these new curvatures are a promising direction for further research in this area. In particular, the associated Ricci flow described in Section \ref{S5: Discrete Ricci Flow} seems to be more amenable to theoretical study than previously considered discrete Ricci flows.  
\\
Apart from the theoretical aspects, the resistance curvature may also prove to play an important role in the context of applications of discrete curvature. The resistance curvature assumes a space between the efficiently calculable but somewhat restrictively local Forman-Ricci curvature and the more integrative/expressive but computationally heavier Ollivier-Ricci curvature -- the resistance curvatures $p$ and $\kappa$ combine local (in its definition) and global (in its dependence on the relative resistance) information on the graph structure and are efficient to approximate. Consequently, the resistance curvature may play a role in the further development of discrete curvature as a tool for network and data analysis.
\\
The precise relation between the introduced resistance curvatures and other notions of discrete curvature is an important unresolved question. We discuss some connections in Section \ref{SS4.2: relation to other discrete curvatures}, describe a formal similarity between $\kappa$ and Ollivier-Ricci curvature in expression \eqref{eq: distance definition link resistance curvature}, and between the equilibrium characterization of $p$ (see Appendix \ref{A: alternative definitions}) and Steinerberger's proposal \cite{steinerberger_curvature_2022}; but these are only partial results and, in particular, say nothing about the relation to the many other notions of discrete curvature not discussed here. However, the classification of different notions of discrete curvature and their relations is an ongoing line of research\footnote{One promising unifying principle of discrete curvature is the so-called semigroup characterization \cite{jost_characterizations_2021}; many important discrete curvatures have such a semigroup characterization, and finding this characterization for the resistance curvature would greatly advance the understanding of its relation to other curvatures.} and detailing the precise position of the resistance curvature in this classification is thus out of scope at this point. Nonetheless, some examples where the curvatures differ (such as the hexagonal lattice, Section \ref{SS4.1: constant curvature graphs}), the fact that resistance curvature depends on the full graph structure compared to local definitions otherwise, the variety of results without an immediate translation to other curvatures -- such as the discrete Ricci flow and the alternative definitions in Appendix \ref{A: alternative definitions} -- and discussions with other researchers, suggest that the resistance curvature is not simply a special case of known notions of discrete curvature.
\\
Finally, we expect that our definitions may be translated from the setting of graphs to more general structures such as simplicial complexes, where (Hodge) Laplacians and effective resistances are also defined \cite{kook_simplicial_2018}. The different definitions in Appendix \ref{A: alternative definitions} provide a number of possible starting points to attempt this translation.
%%%

%%%
%%%%%%%%%%%%%%%%%
%  BIBLIOGRAPHY
%%%%%%%%%%%%%%%%
\bibliographystyle{abbrv}
\bibliography{bibliography.bib}

%%%%%%%%%%%%%%%
%   APPENDIX 
%%%%%%%%%%%%%%%
\appendix
\section*{Appendix}\label{A: appendix general}
\section{Properties of relative resistance and resistance
curvatures}\label{A: properties of rel resistance and res curvature}
\subsection{Relative resistance $c_{ij}\omega_{ij}$}\label{A: properties of relative resistance}
We prove the two statements in Property \ref{property: bounds for relative resistance and Foster's theorem}:\\
\textbf{Proof of bounds for the relative resistance:}
We recall Theorem \ref{th: relative resistance and spanning trees} which says that the relative resistance $c_{ij}\omega_{ij}$ of a link $(i,j)$ is equal to the probability that a random spanning tree contains this link. Since $c_{ij}\omega_{ij}$ is thus a probability, it satisfies $0\leq c_{ij}\omega_{ij}\leq 1$. 
\\
We first refine the \emph{lower bound}. By Theorem \ref{th: resistance is distance} we know that the effective resistance is a distance and thus that $\omega_{ij}>0$ for $i\neq j$. We also know that link weights are positive by definition such that we have $c_{ij}\omega_{ij}>0$.  
\\
Next, we consider when equality occurs for the \emph{upper bound}. If $(i,j)$ is a cut link in $G$, then any subgraph of $G$ without $(i,j)$ is disconnected. Hence, all connected subgraphs (including all spanning trees) of $G$ must contain $(i,j)$ and thus $c_{ij}\omega_{ij}=\Pr[(i,j)\in \mathbf{T}]=1$. Hence, if $(i,j)$ is a cut link then $c_{ij}\omega_{ij}=1$.  Conversely, let $c_{ij}\omega_{ij}=1$ for some link and assume that it is not a cut link. But then if we consider the \emph{connected} graph $G'=G\backslash(i,j)$ from which $(i,j)$ has been removed, we arrive at a contradiction: any spanning tree $T$ of $G'$ is also a spanning tree of $G$ (since $\mathcal{L}_T\subseteq\mathcal{L}_{G'}\subseteq\mathcal{L}_G$ and $G'$ is connected) and $T$ does not contain $(i,j)$ since it is a subgraph of $G'$; however, $c_{ij}\omega_{ij}=1$ implies that $(i,j)$ is in all spanning trees of $G$ which contradicts the existence of such a spanning tree $T$. Consequently, we know that if $\omega_{ij}c_{ij}=1$ then $(i,j)$ must be a cut link. Together, we have thus shown that $\omega_{ij}c_{ij}=1\Leftrightarrow (i,j)$ is a cut link. This completes the proof. \hfill$\square$
\\
\textbf{Proof of Foster's Theorem:}
Assume $G$ is connected. We again use Theorem \ref{th: relative resistance and spanning trees} on the relation between $c_{ij}\omega_{ij}$ and random spanning trees to rewrite the sum over all relative resistances:
\begin{align}
\sum_{(i,j)\in \mathcal{L}} c_{ij}\omega_{ij} 
&= \sum_{(i,j)\in\mathcal{L}}\Pr[(i,j)\in\mathbf{T}] \nonumber\\
&\overset{\text{(a)}}{=} \sum_{(i,j)\in\mathcal{L}}\left(\sum_{T\in\mathcal{T}}\Pr[\mathbf{T}=T]\mathbf{1}_{\lbrace(i,j)\in T\rbrace}\right) \nonumber\\
&= \sum_{T\in\mathcal{T}}\Pr[\mathbf{T}=T]\sum_{(i,j)\in\mathcal{L}}\mathbf{1}_{\lbrace(i,j)\in T\rbrace}\label{eq: in between result}\\
&\overset{\text{(b)}}{=}\sum_{T\in\mathcal{T}}\Pr[\mathbf{T}=T](n-1)\nonumber
\\
&= n-1\nonumber.
\end{align}
In step (a) we make the probability expression explicit and in step (b) we use the fact that any spanning tree of a connected graph has $n-1$ links, independent of the tree. In brief, this derivation shows that $\sum_{i\sim j}c_{ij}\omega_{ij}$ equals the average number of links in a spanning tree, which is constant and equal to $n-1$. For a disconnected graph, this result holds for every component independently and as a result we find that the sum over all links in the graph equals $n-\beta(G)$; this completes the proof.\hfill$\square$
\\
Expression \eqref{eq: in between result} can be used for other sums of relative resistances. For instance, since every cycle must at least miss one link in a spanning tree (which is cycle-free) we find that $\sum_{(i,j)\in\lbrace \text{cycle}\rbrace}\mathbf{1}_{\lbrace (i,j)\in T\rbrace} \leq \vert \lbrace\text{cycle}\rbrace\vert-1$ for any cycle $\lbrace \text{cycle}\rbrace\subseteq\mathcal{L}$. Consequently, the sum of relative resistances over a cycle satisfies
$$
\sum_{(i,j)\in\lbrace \text{cycle}\rbrace}c_{ij}\omega_{ij}\leq \vert \lbrace\text{cycle}\rbrace\vert-1 \text{~for all cycles.}
$$ 
Similarly, since all spanning trees must have at least one link in a given cut -- this is any set of links such that removing the set disconnects the graph\footnote{Often only minimal (with respect to subsets) cuts are considered, which correspond to the links between a bipartition of the graph \cite{bollobas_modern_1998}.} --, we find
\begin{equation}\label{eq: cut bounds}
\sum_{(i,j)\in\lbrace \text{cut}\rbrace}c_{ij}\omega_{ij}\geq 1 \text{~for all cuts.}
\end{equation}
These bounds are an improvement over the straightforward bounds using $0< c_{ij}\omega_{ij}\leq 1$.
%%%%%%%%%

%%%%%%%%%
\subsection{Node resistance curvature $p$}\label{A: properties of node resistance curvature p}
We prove the two statements in Property \ref{property: properties for p_i}:
\\
\textbf{Proof of $\sum p_i=\beta$:} Invoking Foster's Theorem, we can write the sum over all node resistance curvatures as
$$
\sum_{i\in\mathcal{N}} p_i = n - \frac{1}{2}\sum_{i\in\mathcal{N}}\sum_{j\sim i}c_{ij}\omega_{ij} = n-\sum_{(i,j)\in\mathcal{L}}c_{ij}\omega_{ij}=\beta,
$$
which completes the proof. \hfill$\square$
\\
To prove the node resistance curvature bounds in Property \ref{property: properties for p_i}, we first show a stronger result:
\begin{property}\label{property: stronger properties for p_i}
The resistance curvature of a node with $d_i\geq 1$ is bounded by
$$
1-\frac{d_i}{2} \leq p_i \leq 1-\frac{\beta(G_i\backslash\lbrace i\rbrace)}{2} \textup{~for any node $i\in\mathcal{N}$},
$$
where $G_i\backslash\lbrace i\rbrace$ is the component $G_i$ with node $i$ removed. Equality is achieved if and only if all links connected to $i$ are cut links, in which case both bounds are equal.
\end{property}
\textbf{Proof:} The \emph{lower bound} follows immediately from the definition \eqref{eq: definition resistance curvature} of resistance curvature and the relative resistance bound $c_{ij}\omega_{ij}\leq 1$ with equality if and only if all incident links are cut links. The \emph{upper bound} follows by summing the relative resistances over cuts: if removing $i$ (which is possible since $d_i\geq 1$) disconnects $G_i$ into $\beta':=\beta(G_i\backslash\lbrace i\rbrace)$ components, this means that the links incident to node $i$ can be partitioned into $\beta'$ sets of links $\lbrace C_k\rbrace_{k=1}^{\beta'}$ which are cuts of $G_i$. By the bound \eqref{eq: cut bounds} for the sum of relative resistances over a cut, we then find
$$
p_i=1-\frac{1}{2}\sum_{k=1}^{\beta'}\sum_{(i,j)\in C_k}c_{ij}\omega_{ij}\overset{\text{(cut bound)}}{\leq}1-\frac{\beta(G_i\backslash\lbrace i\rbrace)}{2}.
$$
Next, we consider when equality occurs in the upper bound. First, if all links connected to $i$ are cut links then we know (a) that their relative resistances are equal to $1$ and (b) that $\beta(G_i\backslash\lbrace i\rbrace)=d_i$ such that $p_i\leq 1-d_i/2$. Since this is equal to the lower bound, we know that equality must hold. To prove \emph{conversely} that equality in the upper bound implies that all links connected to $i$ are cut links, we will make use of basic facts about spanning trees and cycles in graphs, see for instance \cite{bollobas_modern_1998}. Since the upper bound for $p_i$ was based on the cut bound \eqref{eq: cut bounds}, we know that if this upper bound holds with equality as $p_i=1-\beta(G_i\backslash\lbrace i\rbrace)/2$ then the cut bound \eqref{eq: cut bounds} must also hold with equality. By equation \eqref{eq: in between result} from which the cut bound follows, this implies that $\sum_{(i,j)\in C_k}\mathbf{1}\lbrace (i,j)\in T\rbrace=1$ for all spanning trees $T$ and all cuts $C_k$; in other words, we find that \emph{every spanning tree has exactly one link in every cut $C_k$}. If $n\leq 2$ it then easily follows that the link is a cut link, so assume $n>2$. Let $T$ be a spanning tree, $C_k$ a cut and $\ell=T\cap C_k$ the unique link of $T$ in the cut $C_k$. If we assume that $\vert C_k\vert>1$ then we may take another link in the cut $\ell'\in C_k$ and consider the graph $T+\ell'$, i.e. the tree with this new link added. This graph will have a unique cycle $H$ (see \cite{bollobas_modern_1998}) and, since cycles have length at least $3$, we can take a link $\ell''\in H$ which is different from $\ell,\ell'$. But then the graph $T+\ell'-\ell''$ is again a spanning tree of $G$ (see \cite{bollobas_modern_1998}) which contains both $\ell$ and $\ell'$ i.e. two links in the cut $C_k$. Since this is in contradiction with every spanning tree having exactly one link in the cut $C_k$, we know that $\vert C_k\vert=1$ which means $C_k$ is a cut link; this holds for all cuts $k=1,\dots,\beta'$. In other words, if the upper bound holds with equality we know that all links connected to $i$ are cut links. This completes the proof.\hfill$\square$
\\
\textbf{Proof of node resistance curvature bounds:} The upper bound in Property \ref{property: properties for p_i} follow as a corollary of Property \ref{property: stronger properties for p_i} since $\beta(G_i\backslash\lbrace i\rbrace)\geq \beta(G_i)=1$.\hfill$\square$
\\
Property \ref{property: stronger properties for p_i} shows that the resistance curvature is constrained by the combinatorial structure of the graph as neither $d_i$ nor $\beta(G_i\backslash\lbrace i\rbrace)$ depend on the weights $c$ in the graph. Furthermore, this gives an example of nodes which have a nonpositive curvature:
\begin{property}[cut nodes]\label{prop: cut node is negative}
Cut nodes have nonpositive curvature, with zero curvature if and only if their degree is two.
\end{property}
\textbf{Proof:} If $i$ is a cut node, then $\beta(G_i\backslash\lbrace i\rbrace)\geq 2$ and thus by Property \ref{property: stronger properties for p_i} it follows that $p_i\leq 0$. Next, we consider when zero curvature occurs. Any cut node must have $d_i>1$ since otherwise $i$ is disconnected or a pendant/leaf node and thus not a cut node. If $i$ is a cut node with $d_i=2$, then the two links incident on $i$ must be cut links with unit relative resistance and thus $p_i=0$. Conversely, if $p_i=0$ for a cut node, then from Property \ref{property: stronger properties for p_i} we know that $\beta(G_i\backslash\lbrace i\rbrace)=2$ must hold \emph{and} that all links connected to $i$ must be cut links. This implies that $d_i=2$ as required.\hfill$\square$
\\
A related result about cut nodes and the sign of curvature was shown by Fiedler in \cite{fiedler_matrices_2011}, and follows as a corollary of Property \ref{prop: cut node is negative}:
\begin{corollary}
If $G$ has a cut node then either $G$ is a path graph or $G$ has a node with negative node resistance curvature.
\end{corollary}
\textbf{Proof:} Let $G$ be a graph with cut node $i$. Either $i$ has degree $d_i>2$ and thus negative curvature by Property \ref{prop: cut node is negative}, or degree $d_i=2$ and zero curvature. If $d_i=2$ we furthermore know that the incident links are cut links; hence, if we consider a neighbour of $i$, then either (a) this neighbour has degree $d=1$, in which case it is a leaf nodes with $p=1/2>0$, or (b) this neighbour has degree $d=2$ and both incident links are cut links, such that $p=0$ and we may consider its neighbour recursively, or (c) this neighbour has degree $d>2$ and either (c1) all incident links are cut links such that $p=1-d/2<0$ or (c2) not all incident links are cut links, but then the upper bound in Property \ref{property: stronger properties for p_i} is strict and $p<0$. Thus, either we encounter a node with negative curvature or we can recursively consider further neighbours of degree $d=2$ until we reach a node of degree $d=1$ and stop the process, which would correspond to a path graph. This completes the proof.\hfill$\square$
\\~\\
We mention one additional identity for the node resistance curvature which will be used in the rest of this Appendix:
The product of the Laplacian and resistance matrix of a \emph{connected} graph equals
\begin{equation}\label{eq: QOmega}
Q\Omega = -2I + 2\mathbf{p}\mathbf{u}^T,
\end{equation}
where $\mathbf{u}=(1,\dots,1)^T$ is the all-one vector and $\Omega$ the resistance matrix. This identity follows from \emph{Fiedler's identity} \cite{devriendt_effective_2020,fiedler_matrices_2011} and the work of Bapat \cite{bapat_resistance_2004} which relate the Laplacian and resistance matrices of a graph. In particular, we note that this implies that $\diag(Q\Omega) = 2(\mathbf{p}-\mathbf{u})$.
%%%%%%

%%%%%%
\subsection{Link resistance curvature $\kappa$}\label{A: properties of link resistance curvature kappa}
\textbf{Proof of Property \ref{property: bounds for kappa_ij}:} The bounds for the link resistance curvature $\kappa_{ij}$ in Property \ref{property: bounds for kappa_ij} follow immediately from the bounds for the node resistance curvature $p_i$ in Property \ref{property: properties for p_i}. \hfill$\square$
\\
We can also prove the following stronger result:
\begin{property}[link resistance curvature bounds]\label{property: stronger bounds for kappa_ij}
The link resistance curvature is bounded by
$$
\begin{cases}
\kappa_{ij}\geq \tfrac{1}{\omega_{ij}}(4-d_i-d_j)\\
\kappa_{ij}\leq \tfrac{1}{\omega_{ij}}\left(6-2\beta(G_i\backslash(i,j))-\beta(G_i\backslash\lbrace i,j\rbrace)\right)
\end{cases}
$$
with both bounds achieved (and equal) if and only if all links incident to $i$ and $j$ are cut links, and where $G_i\backslash(i,j)$ is the graph with link $(i,j)$ removed and $G_i\backslash\lbrace i,j\rbrace$ the graph with nodes $i,j$ removed.
\end{property}
\textbf{Proof:} The \emph{lower bound} for the link resistance curvature follows directly from the lower bound $p_i\geq 1-d_i/2$ on the node resistance curvature in Property \ref{property: properties for p_i} and the bound $c_{ij}\omega_{ij}\leq 1$ in Property \ref{property: bounds for relative resistance and Foster's theorem}, with equality only if all incident links are cut links.
\\
We continue with the \emph{upper bound}. For simplicity, we further assume that $G=G_i$ is a connected graph -- this is without loss of generality since any link will be contained in a unique connected component. We can write the link resistance curvature as
\begin{equation}\label{eq: rewritten link curvature}
\frac{2(p_i+p_j)}{\omega_{ij}} = \frac{1}{\omega_{ij}}\left(4 - \sum_{l\in\mathcal{L}_{ij}}c_l\omega_l - 2c_{ij}\omega_{ij} \right),
\end{equation}
where $\mathcal{L}_{ij}:=\lbrace (x,k)\in\mathcal{L}\backslash(i,j)\mid x\in\lbrace i,j\rbrace\rbrace$ are the links incident on $(i,j)$. We note that removing nodes $i$ and $j$ from $G$ is equal to removing the links $\mathcal{L}_{ij}$ from $G$ and then removing the nodes $i,j$ from the resulting graph. Hence, every connected component in $G\backslash\lbrace i,j\rbrace$ is disconnected from $G$ by some cut $C_k\subseteq \mathcal{L}_{ij}$, which gives a partition of $\mathcal{L}_{ij}$ into cuts. For each cut, we may then invoke the cut bound \eqref{eq: cut bounds} which yields
$$
\sum_{l\in\mathcal{L}_{ij}}c_l\omega_l = \sum_{k=1}^{\beta(G\backslash\lbrace i,j\rbrace)}\sum_{l\in C_k}c_l\omega_l \overset{\text{(cut bound)}}{\geq} \beta(G\backslash\lbrace i,j\rbrace),
$$
with equality if and only if all cuts consist of single links (see proof of Property \ref{property: stronger properties for p_i}) and thus if all links in $\mathcal{L}_{ij}$ -- and thus also $(i,j)$ -- are cut links. Second, we know that the relative resistance satisfies $c_{ij}\omega_{ij}\geq 0$ and that $c_{ij}\omega_{ij}=1$ if and only if $(i,j)$ is a bridge link. Hence in general we have $c_{ij}\omega_{ij}\geq [\beta(G\backslash(i,j))-1]$ where $G\backslash(i,j)$ has the link $(i,j)$ removed; again, equality only holds if this is a cut link. Introducing the bounds for these link terms into \eqref{eq: rewritten link curvature} yields the proposed upper bound and conditions for equality, and thus completes the proof.\hfill$\square$
%%%%

%%%%
\section{Alternative definitions for the resistance curvature}\label{A: alternative definitions}
Definition \ref{eq: definition resistance curvature} introduces the node resistance curvature $p$ based on the sum of the relative resistances of links incident on a node. However, this is just one of several equivalent ways to define the resistance curvature, as noted in \cite[Appendix I]{devriendt_variance_2021}. Since these alternative definitions highlight different aspects of the resistance curvature and illustrate how the study of $p$ can benefit from the rich theory of effective resistances, we repeat the definitions and their context from \cite{devriendt_variance_2021}. Each definition also naturally comes with a new expression for the value $\sigma^2=\tfrac{1}{2}\mathbf{p}^T\Omega\mathbf{p}$ introduced in \cite{devriendt_variance_2021} and also discussed in Appendix \ref{A: res Ricci flow}. At the end of this section, we prove Proposition \ref{proposition: distance definitions for resistance curvature} from Section \ref{SS: average distance characterization}.
\\
Many results in this section follow from the identity \eqref{eq: QOmega} in Appendix \ref{A: properties of rel resistance and res curvature} which, we recall, is implied by Fiedler's identity \cite{devriendt_effective_2020}. We repeat the identity for ease of presentation:
\begin{equation}
    Q\Omega = -2I + 2 \mathbf{u}\mathbf{p}^T \tag{Eq. \eqref{eq: QOmega}},
\end{equation}
and recall that $\mathbf{u}=(1,\dots,1)^T$ is the all-one vector and $\Omega$ the resistance matrix. In the rest of this section, the graphs are assumed to be connected.
\subsection{Distance characterization}
A first family of definitions for the node resistance curvature follows immediately from identity \eqref{eq: QOmega}. From the diagonal entries $(Q\Omega)_{ii}$ we retrieve definition \eqref{eq: definition resistance curvature} for the node resistance curvature as
$$
p_i = 1-\frac{1}{2}\sum_{j\sim i}c_{ij}\omega_{ij} \quad\text{~summarized as~}\quad \mathbf{p}=\frac{1}{2}\diag(Q\Omega)+\mathbf{u}.
$$
Looking at the off-diagonal entry $(Q\Omega)_{ix}$ instead, we find 
\begin{equation}\label{eq: distance difference definition}
p_i = \frac{k_i}{2}\left(\omega_{ix} - \sum_{j\sim i}\frac{c_{ij}}{k_i}\omega_{jx}\right) \text{~for any $x\neq i$}.
\end{equation}
Definition \eqref{eq: distance difference definition} shows that the node resistance curvature of $i$ is related to the difference in (resistance) distance from $i$ to $x$ and the average distance from its neighbours to $x$; furthermore, this definition is independent of $x$. Combining definition \eqref{eq: definition resistance curvature} and expression \eqref{eq: distance difference definition} with $x$ ranging over the neighbours of $i$ with probability $c_{ij}/k_i$, we obtain the expression
$$
p_i = \frac{1}{2} - \frac{1}{4}\sum_{k,j\sim i}\frac{c_{ij}c_{ik}}{k_i}\omega_{kj}.
$$
In other words, the resistance curvature of a node is related to the average resistance distance between pairs of neighbours.
\\
In a similar way, if we left-multiply identity \eqref{eq: QOmega} by $\Omega$ and use the fact that $\Omega\mathbf{p}=2\sigma^2$ as follows from \eqref{eq: equilibrium definition for p} in the next section, we find
$$
\sigma^2 = \frac{1}{4}\sum_{j\sim i}c_{ij}(\omega_{ix}-\omega_{jx})^2 \text{~for any $x\neq i$}.
$$
This expression shows that $\sigma^2$ is related to the (squared) difference in distance from a node $x$ to the two end points of each link in a graph, independent of $x$.
\subsection{Equilibrium characterization}
Starting again from identity \eqref{eq: QOmega} and right-multiplying with the resistance curvature vector $\mathbf{p}$, we find that $Q\Omega\mathbf{p} = 0 \Rightarrow \Omega\mathbf{p}\in\ker(Q)$. Since $\ker(Q)=\spn(\mathbf{u})$ for a connected graph \cite{chung_spectral_1997, devriendt_variance_2021} this means that $\Omega\mathbf{p}$ is a constant vector and, including the unit-sum property\footnote{Without the unit-sum property, this definition is underdetermined.}, that $\Omega\mathbf{p}=2\sigma^2\mathbf{u}$. Since the resistance matrix is invertible \cite{devriendt_effective_2020}, this can be written as
\begin{equation}\label{eq: equilibrium definition for p}
\mathbf{p} = \frac{\Omega^{-1}\mathbf{u}}{\mathbf{u}^T\Omega^{-1}\mathbf{u}} \quad\text{~and~}\quad\sigma^2 = \frac{1}{2\mathbf{u}^T\Omega^{-1}\mathbf{u}}.
\end{equation}
As noted in Section \ref{SS4.2: relation to other discrete curvatures}, Steinerberger \cite{steinerberger_curvature_2022} defines a discrete node curvature $\mathbf{w}$ on graphs based on the equation $D\mathbf{w}=\mathbf{u}$ for the shortest-path distance matrix $D$. Taking the resistance distance instead, this would correspond to $\mathbf{w}=\mathbf{p}/(2\sigma^2)$.
\\
We note that combining \eqref{eq: equilibrium definition for p} and identity \eqref{eq: QOmega} for $Q\Omega$ retrieves Bapat's expression for the inverse resistance matrix of a graph \cite{bapat_resistance_2004}:
$$
\Omega^{-1} = -\frac{1}{2}Q + \frac{\mathbf{p}\mathbf{p}^T}{2\sigma^2}.
$$
\subsection{Geometric characterization}
The effective resistance is a \emph{negative type metric} which means there exists an isometric embedding $\varphi:\mathcal{N}\rightarrow\mathbb{R}^{n-1}$ of the square root effective resistance, i.e. with $\Vert\varphi(i)-\varphi(j)\Vert^2 = \omega_{ij}$ for all $i\neq j$. Furthermore, as discovered by Fiedler \cite{fiedler_matrices_2011} (see also \cite{devriendt_effective_2020}) this embedding maps the nodes of a graph to the vertices of a \emph{hyperacute simplex $S$ in $\mathbb{R}^{n-1}$} -- a simplex with nonobtuse (acute or right) interior dihedral angles. This embedding can be linearly extended to functions on the nodes as 
$$
\varphi(f) := \sum_{i\in\mathcal{N}}f(i)\varphi(i) \text{~for any $f:\mathcal{N}\rightarrow\mathbb{R}$},
$$
where a function $f$ is thus the coordinate of a point $\varphi(f)$ in $\mathbb{R}^{n-1}$. As shown by Fiedler \cite{fiedler_matrices_2011} (see also \cite[Appendix I]{devriendt_variance_2021}), we find that the node resistance curvature and $\sigma^2$ satisfy
$$
\Vert\varphi(p)-\varphi(i)\Vert^2 = \sigma^2 \text{~for all $i\in\mathcal{N}$}.
$$
Geometrically, this relates to the circumsphere of $S$ -- the unique sphere passing through all vertices of $S$ -- as follows: \emph{the node resistance curvature is the (unit-sum) coordinate of the center of the circumsphere of $S$ and $\sigma$ is its radius}.
\subsection{Variational characterization}
In \cite{devriendt_variance_2021} the following \emph{variational characterization} of $p$ as the solution to an optimization problem is shown:
$$
\mathbf{p} = \underset{\mathbf{f}:\mathbf{u}^T\mathbf{f}=1}{\argmax}~ \frac{1}{2}\mathbf{f}^T\Omega\mathbf{f} \text{,~~with optimal value $\sigma^2$}.
$$
In \cite{devriendt_variance_2021} we furthermore show that the `variance' of a distribution $\mathbf{f}$ (a nonnegative unit-sum vector) on a graph can be measured by the quadratic product $\var(\mathbf{f}):=\tfrac{1}{2}\mathbf{f}^T\Omega\mathbf{f}$. This leads to expression \eqref{eq: curvature as variance} of the node resistance curvature in terms of the variance of distribution $\mathbf{f}=\brho_{t,i}$.
\subsection{Average distance characterization}\label{AA: proof of average distance characterization}
We prove Proposition \ref{proposition: distance definitions for resistance curvature} which expresses the node and link resistance curvature in terms of the average distance between diffusion distributions around the node or link.\\
\textbf{Proof of Proposition \ref{proposition: distance definitions for resistance curvature}}
Introducing the definitions of the average distance between random nodes \eqref{eq: distance difference definition} and the neighbourhood diffusion distribution $\brho_{t,i}=\exp(-Qt)\mathbf{e}_i$ in the righthandside of expression \eqref{eq: distance definition node resistance curvature}, we obtain
$$
\lim_{t\rightarrow 0}\left(1 - \frac{1}{4t}\mathbb{E}(\omega_{N_tM_t})\right) = \lim_{t\rightarrow 0} \left(1 - \frac{1}{4t}\mathbf{e}_i^T\exp(-Qt)\Omega\exp(-Qt)\mathbf{e}_i\right) \text{~for $N_t,M_t\sim\brho_{t,i}$}.
$$
By definition of the matrix exponential as the power series $\exp(A)=\sum_{k=0}^\infty A^k /(k!)$, the leading order term of the neighbourhood distribution in the $t\rightarrow 0$ limit is $(I-Qt)\mathbf{e}_i$. Making use of identity \eqref{eq: QOmega} for $Q\Omega$, we then find
$$
(I-Qt)\Omega(I-Qt) = \Omega +4I - (2\mathbf{p}\mathbf{u}^T+2\mathbf{u}\mathbf{p}^T)t -2Qt^2,
$$
and thus
$$
\lim_{t\rightarrow 0}\left(1-\frac{1}{4t}\mathbb{E}(\omega_{N_tM_t})\right) = \lim_{t\rightarrow 0} \left(p_i + \frac{1}{2}k_it\right) = p_i \text{~for $N_t,M_t\sim\brho_{t,i}$},
$$
as required. Similarly, for the link resistance curvature we find that the righthandside in expression \eqref{eq: distance definition link resistance curvature} is equal to
$$
\lim_{t\rightarrow 0}\frac{1}{t}\left(1-\frac{\mathbb{E}(\omega_{N_tM_t})}{\omega_{ij}}\right) = \lim_{t\rightarrow 0}\left(\frac{2(p_i+p_j)}{\omega_{ij}} - \frac{2c_{ij}t}{\omega_{ij}}\right) = \kappa_{ij} \text{~for $N_t\sim\brho_{t,i}$ and $M_t\sim\brho_{t,j}$}
$$
as required. This completes the proof.\hfill$\square$
%%%%

%%%%
\section{Resistance curvature and transitivity}\label{A: res curvature and transitivity}
Certain graph symmetries have strong implications for the resistance curvatures in a graph. We show this in particular for node and link transitivity\footnote{These are usually called vertex and edge transitivity.}.
\\
A permutation $\pi:\mathcal{N}\rightarrow\mathcal{N}$ of the nodes of a graph is called an \emph{automorphism} if it preserves the graph structure, i.e. $\pi(i)\sim \pi(j)\Leftrightarrow i\sim j$. A graph is called \emph{node transitive}\footnote{This name refers to the fact that the automorphism group, with automorphisms as group elements and composition `$\circ$' as group operation, acts transitively on the node set.} if for every two nodes $i,j$ there exists an automorphism with $\pi(i)=j$; intuitively, this means that all nodes are indistinguishable in the graph since any two nodes may be interchanged without changing the graph structure. Similarly, a graph is called \emph{link transitive} if for every pair of links $(i,j)$ and $(x,y)$ there exists an automorphism with $(\pi(i),\pi(j))=(x,y)$ (as unordered tuples). See for instance \cite{biggs_algebraic_1974} for properties of such graphs. Some common examples of graphs which are node and link transitive are cycle graphs, the complete graph, the hypercube graph and the Platonic graphs (graph skeletons of the Platonic solids); furthermore all \emph{Cayley graphs} are node transitive \cite{biggs_algebraic_1974, godsil_algebraic_2001}. For graphs with node and/or link transitivity, we can find the resistance curvature exactly as follows:
\begin{property}[transitivity and curvature]\label{prop: transitive graphs}
Let $G$ be a finite connected graph on $n$ nodes and $m$ links and constant link weights $c$, then the following hold:
\begin{itemize}
    \item if $G$ is node transitive, then it has constant resistance curvature equal to $p=1/n$,
    \item if $G$ is link transitive with nodes of (possibly equal) degree $r_1$ and $r_2$, then it has constant link resistance curvature equal to $\kappa=\tfrac{4cm}{n-1} - c(r_1+r_2)$, and
    \item if $G$ is node and link transitive, then it has constant node and link resistance curvatures equal to $p=1/n$ and $\kappa=2c\rho$ where $\rho=m/{{n}\choose{2}}$ is the link density.
\end{itemize}
\end{property}

\textbf{Proof:} 
\textbf{Node transitive} By \cite[Proposition 15.2]{biggs_algebraic_1974}, a permutation $\pi$ of the nodes of a graph is an automorphism if and only if its corresponding permutation matrix $P$ leaves the Laplacian matrix invariant, as $PQP^T=Q$. Consequently, the structure-preserving row and column permutations of the Laplacian of a node transitive graph will act transitively on the row and column set and we say that the Laplacian is row and column transitive. Since $G$ is \emph{connected}, the pseudoinverse Laplacian $Q^\dagger$ is the inverse of the Laplacian in the space $\spn(\mathbf{u})^{\perp}$ with the permutation-invariant constant vector $\mathbf{u}=(1,\dots,1)^T$ and we find that $PQP^T=Q\Rightarrow PQ^\dagger P^T=Q^\dagger$ from:
\begin{align*}
Q^\dagger &= \left(Q+\frac{\mathbf{u}\mathbf{u}^T}{n}\right)^{-1} - \frac{\mathbf{u}\mathbf{u}^T}{n} \quad\quad\text{(see also \cite{ghosh_minimizing_2008})}
\\
&= \left(PQP^T + \frac{\mathbf{u}\mathbf{u}^T}{n}\right)^{-1} - \frac{\mathbf{u}\mathbf{u}^T}{n} \text{~for all $P$ with $PQP^T=Q$}
\\
&= P^T\left[\left(Q+\frac{\mathbf{u}\mathbf{u}^T}{n}\right)^{-1} - \frac{\mathbf{u}\mathbf{u}^T}{n} \right]P \quad\quad\text{~(since $P\mathbf{u}=\mathbf{u}$ and $P^{-1}=P^T$)}
\\
&= P^TQ^\dagger P\text{~for all $P$ with $PQP^T=Q$.}
\end{align*}
This also implies that the pseudoinverse Laplacian $Q^\dagger$ is row and column transitive and consequently that $Q^\dagger$ has a constant diagonal: for every two nodes $i,j$ there exists a permutation matrix $P$ such that $PQ^\dagger P^T=Q^\dagger$ and $P\mathbf{e}_i=\mathbf{e}_j$, and thus 
$$
(Q^\dagger)_{ii} = \mathbf{e}_i^T Q^\dagger \mathbf{e}_i = \mathbf{e}_i^T PQ^\dagger P^T \mathbf{e}_i = \mathbf{e}_j^TQ^\dagger \mathbf{e}_j = (Q^\dagger)_{jj} \text{~for all $i,j$.}
$$
Thus, the diagonal $\bzeta = \diag(Q^\dagger)$ satisfies $P\bzeta = \bzeta$ for any permutation $P$. Since we can write the resistance matrix as $\Omega = \mathbf{u}\bzeta^T + \bzeta\mathbf{u}^T-2Q^\dagger$ following definition \eqref{eq: definition effective resistance} (see also \cite{devriendt_effective_2020}), this implies that $PQP^T=Q\Rightarrow P\Omega P^T=\Omega$ and that the resistance matrix is also row and column transitive.
\\
From the identity $Q\Omega=-2I+2\mathbf{p}\mathbf{u}^T$ in expression \eqref{eq: QOmega}, we then find that
$$
\mathbf{p}=\frac{1}{2}\diag(Q\Omega)+\mathbf{u}=\frac{1}{2}\diag(PQP^TP\Omega P^T)+\mathbf{u}=P\mathbf{p}
$$
for all structure-preserving row/column permutations $P$ of $Q$. Since these permutations act transitively on the rows, $p$ must be constant. By Property \ref{property: properties for p_i} we then know that $\sum_i p_i=np=1$ and thus that $p_i=1/n$ for all nodes.
\\
\textbf{Link transitive but not node transitive} If a graph is link transitive and not node transitive, then it must be bipartite on $\mathcal{V}_1,\mathcal{V}_2$ and the automorphism group acts transitively on these partitions, i.e. the nodes are indistinguishable inside the partitions, see \cite[Lemma 3.2.1]{godsil_algebraic_2001}. In particular, this implies that all ($n_1$) nodes in $\mathcal{V}_1$ have the same degree $r_1$ and all ($n_2$) nodes in $\mathcal{V}_2$ have degree $r_2$. The number of links is then equal to $m=n_1r_1=n_2r_2$. By Foster’s theorem and link transitivity -- such that $c_{ij}\omega_{ij}$ must be equal for all links -- we find that the effective resistance of every is equal and given by $\omega = \tfrac{n-1}{n_1r_1c} = \tfrac{n-1}{n_2r_2c}$. The node curvatures in the two sets then follow as $p_1=1-\tfrac{n-1}{2n_1}$ and $p_2=1-\tfrac{n-1}{2n_2}$ and the link resistance curvature is then calculated from their sum resulting in the proposed formula for $\kappa$. 
\\
\textbf{Link and node transitive} If the graph is link transitive and node transitive, then $r_1=r_2$ in the previous derivation, and from $c(r_1+r_2)=4mc/n$ we then find the proposed link resistance curvature $2rc/(n-1)$. An alternative derivation uses the fact that if the graph is node transitive, then $p_i=1/n$ as before and that Foster’s Theorem together with link transitivity implies that $\omega=\tfrac{2(n-1)}{nrc}$ such that $\kappa=2rc/(n-1)$ as required. This completes the proof.\hfill$\square$
%%%%

%%%%
\section{Correspondences between resistance curvature and established discrete curvatures}\label{A: correspondence between res curvature and other curvatures}
\subsection{Combinatorial curvature}\label{A: correspondence with combinatorial curvature}
We show formula \eqref{eq: resistance curvature is expected combinatorial curvature} which relates the node resistance curvature and combinatorial curvature. Starting from Definition \ref{eq: definition resistance curvature} for the resistance curvature and using Theorem \ref{th: relative resistance and spanning trees} for the relation between relative resistances and spanning trees, we find
\begin{align*}
p_i &= 1 - \frac{1}{2}\sum_{j\sim i}c_{ij}\omega_{ij}\\
&= 1 - \frac{1}{2}\sum_{j\sim i}\sum_{T\in\mathcal{T}}\Pr[\mathbf{T}=T]\mathbf{1}_{\lbrace (i,j)\in T\rbrace}\quad\text{~(by Theorem \ref{th: relative resistance and spanning trees})}\\
&= \sum_{T\in\mathcal{T}}\Pr[\mathbf{T}=T]\left(1 - \frac{1}{2}\sum_{j\sim i}\mathbf{1}_{\lbrace (i,j)\in T\rbrace}\right)\\
&= \sum_{T\in\mathcal{T}}\Pr[\mathbf{T}=T] \left(1 - \frac{d_i^{(T)}}{2}\right) = \mathbb{E}[p_i^{(co)}(\mathbf{T})],
\end{align*}
where $d_i^{(T)}$ is the degree of node $i$ in spanning tree $T$. This proves formula \eqref{eq: resistance curvature is expected combinatorial curvature}.
\subsection{Ollivier-Ricci curvature}\label{A: correspondence with Ollivier curvature}
We prove Proposition \ref{prop: bounds for OR curvature} which relates the link resistance curvature to Ollivier-Ricci curvature measured with respect to the effective resistance metric and balls $\mu_{t,i}$ determined by a lazy random walk, i.e. given by the vector $\bmu_{t,i}=(I-Qt)\mathbf{e}_i$. We will assume $G$ to be connected, as both $\kappa$ and $\kappa^{(OR)}$ are determined in each component independently.
\\
\textbf{Proof of Proposition \ref{prop: bounds for OR curvature}:} We start by bounding the Wasserstein distance between the balls $\mu_{t,i}$ and $\mu_{t,j}$:
\begin{equation}\label{eq: bound for wasserstein}
W_1(\mu_{t,i},\mu_{t,j}) = \min_P \tr(P\Omega) \overset{(a)}{\leq} \tr(\bmu_{t,i}\bmu_{t,j}^T \Omega) \overset{(b)}{=} \bmu_{t,i}^T\Omega\bmu_{t,j},
\end{equation}
where for the inequality (a) we have used that $P = \bmu_{t,i}\bmu_{t,j}^T$ is a valid matrix for the Wasserstein distance, i.e. with nonnegative entries and the correct marginals, and where equality (b) invokes properties of the trace operator. Introducing the definition of the lazy random walk for $\bmu_{t,i}$ and $\bmu_{t,j}$, we obtain
$$
W_1(\mu_{t,i},\mu_{t,j}) \leq \mathbf{e}_i^T(I-Qt)\Omega(I-Qt)\mathbf{e}_j \overset{(a)}{=} \omega_{ij} - 2t(p_i+p_j) + 2t^2c_{ij},
$$
where in (a) we use the identity $Q\Omega = -2I+2\mathbf{p}\mathbf{u}^T$ as in \eqref{eq: QOmega}. Introduced into the definition of Ollivier-Ricci curvature with $\omega$ as distance, we then find
\begin{equation}\label{eq: inequality between OR and link resistance curvature}
\kappa_{ij}^{(OR)} =\lim_{t\rightarrow 0^+}\frac{1}{t}\left( 1-\frac{W_1(\mu_{t,i},\mu_{t,j})}{\omega_{ij}}\right) \geq \frac{2(p_i+p_j)}{\omega_{ij}},
\end{equation}
establishing the bound in Proposition \ref{prop: bounds for OR curvature}.
\\
Next, we show that equality is achieved in the case of cut links. If $(i,j)$ is a cut link, then $i$ is a cut node and thus for any two nodes $x,y$ which are disconnected by removal of $i$, we have the triangle \emph{equality} $\omega_{xy}=\omega_{xi}+\omega_{iy}$, see for instance \cite{klein_resistance_1993}; similarly, $j$ is a cut node. In particular, the triangle equality will hold for effective resistances between the neighbours of $i$ and the neighbours of $j$, in other words between the supports $\mathcal{I}:=\supp(\mu_{t,i})$ and $\mathcal{J}:=\supp(\mu_{t,j})$ of the balls around $i$ and $j$ respectively -- where for the lazy random walk we also have $i\in\mathcal{I}$ and $j\in\mathcal{J}$ and $j\in\mathcal{I},i\in\mathcal{J}$ because $i\sim j$. These effective resistances satisfy $\omega_{xy} = \omega_{xi}+\omega_{ij}+\omega_{jy}$ for all $x\in\mathcal{I}\backslash\lbrace j\rbrace,y\in\mathcal{J}\backslash\lbrace i\rbrace$ and consequently, the block matrix $\Omega_{\mathcal{J}\mathcal{I}}$ has the following low-rank decomposition
$$
\Omega_{\mathcal{J}\mathcal{I}} = \bomega_{j}\mathbf{u}^T + \mathbf{u}\bomega_i^T + \omega_{ij}\left(\mathbf{u}\mathbf{u}^T -2\mathbf{e}_i\mathbf{u}^T -2\mathbf{u}\mathbf{e}_j^T +2\mathbf{e}_i\mathbf{e}_j^T\right)
$$
where $\bomega_i$ is the $\vert\mathcal{I}\vert\times 1$ vector with effective resistances to the neighbours of $i$, as $(\bomega_i)_x = \omega_{ix}$ for $x\in\mathcal{I}$ and similarly for $\bomega_j$, and with $\mathbf{u}$ the all-one vectors of appropriate size. Introducing this into the expression for the Wasserstein distance, we retrieve
\begin{align*}
W_1(\mu_{t,i},\mu_{t,j}) &= \min_{P}\tr(P\Omega) 
\\
&\overset{\text{(a)}}{=} \min_{P}\tr(P_{\mathcal{I}\mathcal{J}}\Omega_{\mathcal{J}\mathcal{I}}) \\
&= \min_{P} \left(\mathbf{u}^TP_{\mathcal{I}\mathcal{J}}\bomega_j + \bomega_i^T P_{\mathcal{I}\mathcal{J}}\mathbf{u} + \omega_{ij}\left[\mathbf{u}^TP_{\mathcal{I}\mathcal{J}}\mathbf{u}-2\mathbf{u}^TP_{\mathcal{I}\mathcal{J}}\mathbf{e}_i - 2\mathbf{e}_j^TP_{\mathcal{I}\mathcal{J}}\mathbf{u} + 2\mathbf{e}_j^TP_{\mathcal{I}\mathcal{J}}\mathbf{e}_i\right]\right)
\\
&\overset{\text{(b)}}{=} \min_{P}\left(\bmu_{t,j}^T(\bomega_j -2\omega_{ij}\mathbf{e}_i) + (\bomega_i^T-2\omega_{ij}\mathbf{e}_j)^T\bmu_{t,i} + \omega_{ij}(1+2P_{ji})\right)
\\
&= \bmu_{t,i}^T\Omega\bmu_{t,j} - 2\omega_{ij} \bmu_{t,j}^T\mathbf{e}_i\mathbf{e}_j^T\bmu_{t,i} + 2\omega_{ij}\min_{P}P_{ji}
\\
&\overset{\text{(c)}}{\geq} \bmu_{t,j}^T\Omega\bmu_{t,i} -2t^2c_{ij}
\end{align*}
where in step (a) we use that a nonnegative matrix $P$ with $\supp(P\mathbf{u})=\mathcal{I}$ and $\supp(P^T\mathbf{u})=\mathcal{J}$ can only have nonzero entries in $\mathcal{I}\times\mathcal{J}$, in step (b) we introduce the marginals of $P$ and in step (c) we use nonnegativity of $P_{ji}$ and the fact that $(\bmu_{t,i})_{j}=(\bmu_{t,j})_i=c_{ij}t$ and $c_{ij}\omega_{ij}=1$ since the link is a cut link. Introducing the definitions of $\bmu_{t,i}$ and $\bmu_{t,j}$, we obtain
$$
W_1(\mu_{t,i},\mu_{t,j}) \geq \mathbf{e}_j^T(I-Qt)\Omega(I-Qt)\mathbf{e}_i-2t^2c_{ij} = \omega_{ij} - 2t(p_i+p_j)
$$
such that
$$
\kappa_{ij}^{(OR)} = \lim_{t\rightarrow 0^{+}}\frac{1}{t}\left(1 - \frac{W_1(\mu_{t,i},\mu_{t,j})}{\omega_{ij}}\right) \leq \frac{2(p_i+p_j)}{\omega_{ij}} \quad\text{~(if $(i,j)$ is a cut link)}.
$$
Combined with the general lower-bound \eqref{eq: inequality between OR and link resistance curvature} for the OR curvature, this proves that equality must hold between $\kappa_{ij}^{(OR)}$ and $\kappa_{ij}$ when $(i,j)$ is a cut link.\hfill$\square$
\subsection{Forman-Ricci curvature}\label{A: correspondence with Forman curvature}
We show Proposition \ref{prop: upper bound for FR curvature}, which relates the link resistance curvature to Forman-Ricci (FR) curvature \eqref{eq: definition FR curvature} with a particular choice of weights.
\\
\textbf{Proof of Proposition \ref{prop: upper bound for FR curvature}:} Starting from the definition of FR curvature \eqref{eq: definition FR curvature} with unit weights $w=1$, we find that
\begin{align*}
\kappa_{ij}^{(FR)} &= 2\left(1-\frac{1}{2}\sum_{k\sim i}1\right) + 2\left(1-\frac{1}{2}\sum_{k\sim j}1\right) = 4 - d_i-d_j \overset{(a)}{\leq} \omega_{ij}\kappa_{ij}
\end{align*} 
where (a) is the lower bound for $\kappa_{ij}$ in Property \ref{property: bounds for kappa_ij}. This completes the proof.\hfill$\square$
%%%%

%%%%
\section{Resistance curvature in Euclidean random graphs}\label{A: res curvature in ERGs} 
\begin{figure}[h!]
    \centering
    \includegraphics[width=0.5\textwidth]{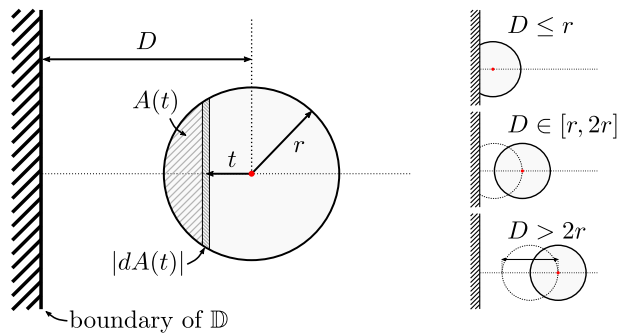}
    \caption{\emph{Illustration of the local neighbourhood around a node in a Euclidean random graph (ERG). By assuming that the boundary curvature is negligible with respect to the curvature of the connection disc around each node as $r\ll r_{\mathbb{D}}$, the boundary can be taken as a straight line. \textbf{\textup{The left figure}} shows a node (red point) at distance $D>r$ from the boundary, the connection disc of radius $r$ around this node, an example of a circle segment at height $t$ with area $A(t)$ and finally $\vert dA(t)\vert$ which is the intersection of points at distance $D-t$ from the boundary and the connection disc. \textbf{\textup{The right figure}} illustrates the three different node regimes determined by their distance $D$ to the boundary.}}\label{fig: setup for ERGs}
\end{figure}

We derive formula \eqref{eq: boundary function formula} as an approximation/model of the expected node resistance curvature in Euclidean random graphs. This formula is based on the heuristic \eqref{eq: heuristic for resistance} for the effective resistance
$$
\hat{\omega}_{ij} = \frac{1}{\lambda S(i)} + \frac{1}{\lambda S(j)},
$$
where $S(i)$ is the area of the intersection of the domain $\mathbb{D}$ and an $r$-radius disc around $i$. In general, this intersection can take many forms depending on the geometry of the boundary and the location of $i$ with respect to the boundary, and to overcome this complexity we will make two further assumptions. First, we assume that the curvature of the boundary (as a $1$D curve in $\mathbb{R}^2$) is negligible with respect to the curvature of the $r$-disc around $i$; if the largest curvature of the boundary is $1/r_{\mathbb{D}}$ then \emph{we assume $r\ll r_{\mathbb{D}}$}. We remark that in the case of a circular domain of radius $R$ as in the main text, we have $r_{\mathbb{D}}=R$ such that the assumption is $r\ll R$. With this assumption, we will be able to consider the boundary locally as a straight line (i.e. when considering the neighbourhood around a node). Second, we assume that every point which lies at distance less than $2r$ from the boundary of $\mathbb{D}$ has a unique closest point on the boundary. This may be formalized using the concept of \emph{reach} \cite{federer_curvature_1959, aamari_estimating_2019} by requiring that the boundary has reach at least $2r$. We remark that for a fixed domain $\mathbb{D}$, these assumptions will automatically be satisfied if we let $r$ be small enough.
\\
As Figure \ref{fig: setup for ERGs} illustrates, with these further assumptions $S(i)$ is determined by the intersection of a disc and a half-plane (due to the straight boundary) and equal to
\begin{equation}\label{eq: area around node}
S(i) = \begin{dcases}
\pi r^2 \text{~if $D_i\geq r$}\\
A(-D_i) \text{~if $0\leq D_i\leq r$}
\end{dcases}\text{~or, in short $S(i)=A(\max\lbrace-r,-D_i\rbrace)$}
\end{equation}
where $A(t)$ is the area of a circular segment\footnote{This is equal to $A(t)=r^2\arccos(t/r)-t\sqrt{r^2-t^2}$, and we find $\vert dA(t)\vert=2\sqrt{r^2-t^2}$.} at height $t$. In other words, the function $S$ only depends on the distance from the boundary, and only variations of node position in the direction towards or away from the boundary can change $S$; we may thus further write $S(D_i)$. We note that $A(-r) = \pi r^2$ and $A(x)+A(-x)=\pi r^2$ for all $x\in[-r,r]$.
\\
Now, we fix a distance $D$ and consider one sample graph $G$ of the ERG on $\mathbb{D}$ to which we add a node $i$ at distance $D_i=D$ from the boundary. For this node, we can write
\begin{align*}
\hat{p}_i &= 1 - \frac{1}{2}\sum_{j\sim i} \hat{\omega}_{ij} \quad\text{(we consider unweighted ERGs)}    
\\
&= 1 - \frac{1}{2}\sum_{j\sim i}\left(\frac{1}{\lambda S(D_i)}+ \frac{1}{\lambda S(D_j)}\right)\quad\text{(introducing \eqref{eq: heuristic for resistance} for $\hat{\omega}$)}
\\
&= 1 - \frac{d_i}{2\lambda S(D)} - \frac{1}{2}\sum_{j\sim i}\frac{1}{\lambda S(D_j)}\quad\text{($D_i=D$ by construction)}
\\
&= 1 - \frac{d_i}{2\lambda S(D)} - \frac{1}{2}\sum_{D'}\frac{\text{\# neighbours of $i$ at distance $D'$ from boundary}}{\lambda S(D')},
\end{align*}
where in the last expression $D'$ sums over $\lbrace D_j: j\sim i\rbrace$. As noted, $S$ only depends on variations in the direction orthogonal to the (straight) boundary. This allows to express the position of the neighbours of $i$ using the difference $D-D_j$, i.e. the distance between $i$ and $j$ measured in the direction orthogonal to the boundary. We can write
\begin{equation}\label{eq: heuristic for p for one sample}
\hat{p}_i = 1 - \frac{d_i}{2\lambda S(D)} - \frac{1}{2}\sum_{t}\frac{\text{\# neighbours of $i$ at distance $D-t$ from the boundary}}{\lambda S(D-t)}
\end{equation}
where we sum over $\lbrace t=D-D_j:j\sim i\rbrace$. We now consider the expectation of $\hat{p}$ over a random realization of the Euclidean random graph in which node $i$ is artificially included at the same point in the domain, at distance $D$ from the boundary. First, the expectation of the degree $d_i$ is proportional to $S(D)$ as a property of the Poisson point process, as
$$
\mathbb{E}(d_i) = \lambda S(D).
$$
Second, as shown in Figure \ref{fig: setup for ERGs} and as a property of the Poisson process, we find that the expected number of neighbours of node $i$ at distance $D-t$ from the boundary is given by $\lambda \vert dA(t)\vert$. Together, we can write the expectation of \eqref{eq: heuristic for p for one sample} as
\begin{equation}\label{eq: heuristic for p in expectation}
\mathbb{E}(\hat{p}_i) = \frac{1}{2}\left(1 - \int_{-r}^{\min(r,D)}\frac{\vert dA(t)\vert}{S(D-t)}\right),
\end{equation}
where the range for $t=D-D_j$ is lower bounded by $-r$, which are the neighbours of $i$ furthest away from the boundary, and upper bounded by $\min(r,D)$, which are the neighbours closest to the boundary. Introducing \eqref{eq: area around node} for the function $S$ retrieves the formula \eqref{eq: boundary function formula} for the expected resistance curvature in Section \ref{SS4.3: resistance curvature for ERGs}. Equation \eqref{eq: heuristic for p in expectation} can be integrated numerically for different values of $D$ as shown in Figure \ref{fig: result for ERGs}.
\\~\\
We can further simplify expression \eqref{eq: heuristic for p in expectation} by distinguishing three different cases for $D$. Most importantly, this shows that nodes which are sufficiently far away from the boundary will have zero curvature $\hat{p}_i=0$ in expectation and thus provides a (heuristic) explanation for the zero resistance curvature tendency in ERGs. The different regimes for $D$ are as follows: 
\\
\textbf{Bulk regime ($D\geq 2r$):} In this case we get $\min(r,D)=r$ and $S(D-t)= \pi r^2$ for all $t\in[-r,r]$, since the neighbours of $i$ are at least at distance $r$ from the boundary. From $\int_{-r}^r\vert dA(t)\vert=A(-r)$ we then find that
$$
\mathbb{E}(\hat{p}_i) = \frac{1}{2}\left(1 - \frac{1}{\pi r^2}\int_{-r}^r \vert dA(t)\vert\right) = 0 \text{~for $D\geq 2r$}.
$$
\textbf{Near-boundary regime ($r\leq D\leq 2r$):} In this case we get $\min(r,D)=r$ and thus neighbours in the range $-r\leq t\leq r$ around $i$. For these neighbours, we find (a) if $-r\leq t \leq D-r$ then $r\leq D-t\leq D+r$ and thus by \eqref{eq: area around node} that $S(D-t) = \pi r^2$ and (b) if $D-r\leq t\leq r$ then $D-r\leq D-t\leq r$ and thus by \eqref{eq: area around node} that $S(D-t) = A(t-D)$. Introduced in \eqref{eq: heuristic for p in expectation} this yields
\begin{align*}
\mathbb{E}(\hat{p}_i) &= \frac{1}{2}\left(1 - \frac{1}{\pi r^2}\int_{-r}^{D-r}\vert dA(t)\vert - \int_{D-r}^r\frac{\vert dA(t)\vert}{A(t-D)}\right)
\\
&= \frac{1}{2}\left(1 - \frac{A(-r)-A(D-r)}{\pi r^2} - \int_{D-r}^r\frac{\vert dA(t)\vert}{A(t-D)}\right)
\\
&= \frac{A(D-r)}{2\pi r^2} - \frac{1}{2}\int_{D-r}^r\frac{\vert dA(t)\vert}{A(t-D)} \text{~for $r\leq D\leq 2r$}.
\end{align*}
\textbf{Boundary regime ($0\leq D\leq r$):} In this case we get $\min(r,D)=D$ and thus neighbours in the range $-r\leq t\leq D$ around $i$. Similar to the near-boundary regime we then find that
\begin{align*}
\mathbb{E}(\hat{p}_i) &= \frac{1}{2}\left(1 - \frac{1}{\pi r^2}\int_{-r}^{D-r}\vert dA(t)\vert - \int_{D-r}^D\frac{\vert dA(t)\vert}{A(t-D)}\right)
\\
&= \frac{A(D-r)}{2\pi r^2} - \frac{1}{2}\int_{D-r}^D\frac{\vert dA(t)\vert}{A(t-D)} \text{~for $0\leq D\leq r$}.
\end{align*}
To summarize, we find the following piecewise expression
$$
\mathbb{E}(\hat{p}(D)) = \begin{dcases}
\frac{A(D-r)}{2\pi r^2} - \frac{1}{2}\int_{D-r}^{D}\frac{\vert dA(t)\vert}{A(t-D)}\text{~if $D\in[0,r]$}\\
\frac{A(D-r)}{2\pi r^2} - \frac{1}{2}\int_{D-r}^{r}\frac{\vert dA(t)\vert}{A(t-D)}\text{~if $D\in[r,2r]$}\\
0 \text{~if $D\geq 2r$}
\end{dcases}
$$
\\
\textit{Remark:} Our analysis of the resistance curvature in ERGs based on $\mathbb{E}(\hat{p})$ not only explains why we may expect zero resistance curvature in the bulk of ERGs, but it also provides an explanation for the `boundary effect' in Figures \ref{fig: constructing ERG and resistance curvature on disc} and \ref{fig: result for ERGs} where we found experimentally that the resistance curvature changes from zero to negative and then positive curvature when nearing the boundary. The analysis above based on $\hat{\omega}$ shows how this boundary effect originates from the different possible local geometries around a node. This analysis could be relevant in other cases, for instance in \cite{willerton_heuristic_2009} where a similar boundary effect was described in the context of the magnitude of graphs (see also \cite{leinster_magnitude_2013}), and \cite{bunch_weighting_2021} where the boundary effect was used for boundary detection and related data analysis tasks.
%%%%

%%%%
\section{Results on the resistance Ricci flow}\label{A: res Ricci flow}
We first show Proposition \ref{proposition: gradient expression for resistance Ricci flow}, which gives a gradient flow expression for the resistance Ricci flow \eqref{eq: definition resistance ricci flow}. Note that we have assumed $G$ to be connected such that the resistance matrix $\Omega$ only contains effective resistances between nodes in the same component.\\
\textbf{Proof of Proposition \ref{proposition: gradient expression for resistance Ricci flow}:} We follow the terminology and approach for matrix differentiation used in \cite{petersen_matrix_2008}. For the derivative of some function $f(\Omega)$ with respect to the resistance matrix $\Omega$, one has to account for the specific structure of this matrix. By symmetry and zero diagonal, we have
$$
\frac{\partial\Omega}{\partial \omega_{ij}} = (1-\mathbf{1}_{\lbrace i=j\rbrace})(\mathbf{e}_i\mathbf{e}_j^T + \mathbf{e}_j\mathbf{e}_i^T).
$$
The derivative of a function $f(\Omega)$ with respect to $\Omega$ is then found from \cite{petersen_matrix_2008}
$$
(\nabla_{\Omega} f)_{ij} = \tr\left(\left[\frac{\partial f}{\partial \Omega}\right]^T\frac{\partial \Omega}{\partial\omega_{ij}}\right),
$$
with matrix gradient $\nabla_{\Omega}$. As the potential function is a trace with partial derivative $\partial\tr(\tfrac{1}{2}\Omega Q\Omega)/\partial\Omega = \frac{1}{2}(Q\Omega+\Omega Q)$, we find
\begin{align*}
(\nabla_{\Omega}\tr(\tfrac{1}{2}\Omega Q\Omega))_{ij} &= \tr\left(\tfrac{1}{2}(1-\mathbf{1}_{\lbrace i=j\rbrace})(Q\Omega + \Omega Q)(\mathbf{e}_i\mathbf{e}_j^T+\mathbf{e}_j\mathbf{e}_i^T)\right)
\\
&= (1-\mathbf{1}_{\lbrace i=j\rbrace})\left(\mathbf{e}_i^TQ\Omega \mathbf{e}_j + \mathbf{e}_j^T Q\Omega \mathbf{e}_i\right)
\\
&= (1-\mathbf{1}_{\lbrace i=j\rbrace})(2p_i+2p_j-4.\mathbf{1}_{\lbrace i=j\rbrace})\quad\text{(by expression \eqref{eq: QOmega} for $Q\Omega$)}
\\
&= \begin{cases}
2(p_i+p_j) \text{~if $i\neq j$}\\
0\text{~otherwise}
\end{cases}
\end{align*}
This confirms that the gradient corresponds to minus the resistance Ricci flow \eqref{eq: definition resistance ricci flow}, which completes the proof of Proposition \ref{proposition: gradient expression for resistance Ricci flow}.\hfill$\square$
\\
\textit{Remark:} with identity \eqref{eq: QOmega} for the product $Q\Omega$, the potential function can also be written as
$$
\tr\left(\tfrac{1}{2}\Omega Q\Omega\right) = \tr\left(\Omega[\mathbf{p}\mathbf{u}^T-I]\right) = \mathbf{u}^T\Omega \mathbf{p}.
$$
As shown in \cite{devriendt_effective_2020} and \cite[Appendix I]{devriendt_variance_2021} and repeated in Appendix \ref{A: alternative definitions}, the node resistance curvature vector $\mathbf{p}$ satisfies the following equations
$$
\begin{dcases}
\Omega\mathbf{p} = 2\sigma^2 \mathbf{u} \text{, where}\\
\sigma^2 = \tfrac{1}{2}\mathbf{p}^T\Omega\mathbf{p}
\end{dcases} \text{~or, equivalently}
\begin{dcases}
\mathbf{p} = \frac{1}{2\sigma^2}\Omega^{-1}\mathbf{u}\text{, where}\\
\sigma^2 = \frac{1}{2\mathbf{u}^T\Omega^{-1}\mathbf{u}}
\end{dcases}
$$
which allows the potential to be written as $\tr(\tfrac{1}{2}\Omega Q\Omega) = 2n\sigma^2$. The expressions for $\mathbf{p}$ and $\sigma^2$ in terms of the inverse resistance matrix are particularly relevant in the context of magnitude \cite{leinster_magnitude_2013}.
\\~\\
Next, we prove Proposition \ref{propos: Laplacian form for Ricci flow} which gives expression \eqref{eq: dQ=2Q^2} for the resistance Ricci flow in terms of Laplacian matrices. We again assume $G$ to be connected for the proof, but the general result follows by combining the Laplacians of each component into a block-diagonal Laplacian matrix.
\\
\textbf{Proof of Proposition \ref{propos: Laplacian form for Ricci flow}:} We will start from expression \eqref{eq: dQ=2Q^2} and show that this flow of Laplacian matrices is equivalent to \eqref{eq: definition resistance ricci flow} as a flow of effective resistances. First, for the derivative of the inverse of a matrix we find that $dA^{-1}/dt = -A^{-1}(dA/dt)A^{-1}$. Writing the Laplacian pseudoinverse as $Q^\dagger = (Q+\mathbf{u}\mathbf{u}^T/n)^{-1}-\mathbf{u}\mathbf{u}^T/n$ (see \cite{ghosh_minimizing_2008}) and considering the flow $dQ/dt=2Q\diag(\mathbf{p})Q$, this yields
\begin{align*}
\frac{dQ^\dagger}{dt} &= -Q^\dagger\frac{dQ}{dt}Q^\dagger
\\
&= -2Q^\dagger Q\diag(\mathbf{p})QQ^\dagger
\\
&= -2\left(I-\frac{\mathbf{u}\mathbf{u}^T}{n}\right)\diag(\mathbf{p})\left(I-\frac{\mathbf{u}\mathbf{u}^T}{n}\right).
\end{align*}
For the change of effective resistances, we then find
$$
\frac{d\omega_{ij}}{dt} = (\mathbf{e}_i-\mathbf{e}_j)^T\frac{dQ^\dagger}{dt}(\mathbf{e}_i-\mathbf{e}_j) = -2(p_i+p_j) \text{~for all $i\neq j$}
$$
which confirms that \eqref{eq: dQ=2Q^2} and \eqref{eq: definition resistance ricci flow} are equivalent expressions for the resistance Ricci flow.\hfill$\square$
\end{document}